\documentclass[11pt,a4paper]{amsart}
\usepackage[margin=2.7cm]{geometry}
\usepackage{amssymb,amsmath,amsthm}
\usepackage{graphicx}
\usepackage{overpic}
\usepackage[table]{xcolor}

\usepackage{chngcntr}
\usepackage{lscape}
\usepackage{rotating}
\usepackage{booktabs}
\usepackage{units}

\usepackage{epstopdf}

\usepackage{multirow}
\usepackage{calc}
\usepackage{subcaption}

\usepackage{tikz}
\usepackage{tikz-cd}
\usepackage{pgfplots}
\usepackage{pgfplotstable}
\usetikzlibrary{3d,calc}
\usetikzlibrary{shapes, spy, scopes, matrix, datavisualization.formats.functions, arrows, datavisualization, plothandlers, plotmarks, svg.path, positioning}
\usetikzlibrary{spy}
\usetikzlibrary{svg.path}
\usetikzlibrary{positioning}
\pgfplotsset{compat=newest}
\pgfplotsset{plot coordinates/math parser=false}
\pgfplotsset{scaled y ticks = false, tick label style={/pgf/number format/fixed}}
\pgfplotsset{scaled x ticks = false, tick label style={/pgf/number format/fixed}}
\usepgfplotslibrary{groupplots}

\usepackage{float}

\usepackage[bookmarksopen, pdfstartview={FitV}, linkbordercolor=white]{hyperref}

\usepackage[]{algorithm2e}

\usepackage{todonotes}

\usepackage[capitalize,nameinlink]{cleveref}
	
\newtheoremstyle{mystyle}
  {}
  {}
  {\itshape}
  {}
  {\bfseries}
  {.}
  { }
  {}%
	
\theoremstyle{mystyle}
\newtheorem{theorem}{Theorem}[section]
\theoremstyle{plain}

\theoremstyle{plain}
\newtheorem{lemma}[theorem]{Lemma}
\theoremstyle{plain}
\newtheorem{remark}[theorem]{Remark}
\theoremstyle{plain}

\theoremstyle{plain}
\newtheorem{corollary}[theorem]{Corollary}
\theoremstyle{plain}
\newtheorem{proposition}[theorem]{Proposition}
\theoremstyle{plain}

\theoremstyle{plain}
\newtheorem{assumption}[theorem]{Assumption}
\theoremstyle{plain}

\theoremstyle{plain}
\newtheorem{conjecture}[theorem]{Conjecture}

\DeclareMathOperator*{\argmin}{\arg\!\min}
\newcommand{\BIGOP}[1]{\mathop{\mathchoice%
{\raise-0.22em\hbox{\huge $#1$}}%
{\raise-0.05em\hbox{\Large $#1$}}{\hbox{\large $#1$}}{#1}}}

\newcommand{\BIGboxplus}{\mathop{\mathchoice%
{\raise-0.35em\hbox{\huge $\boxplus$}}%
{\raise-0.15em\hbox{\Large $\boxplus$}}{\hbox{\large $\boxplus$}}{\boxplus}}}

\crefname{assumption}{Assumption}{Assumptions}
\crefname{notation}{Notation}{Notations}
\crefname{observation}{Observation}{Observations}
\crefname{conjecture}{Conjecture}{Cenjectures}

\def\N{\mathbb N}

\def\R{\mathbb R}

\def\R{\mathbb{R}}
\def\N{\mathbb{N}}

\def\cN{\mathcal{N}}
\def\cO{\mathcal{O}}

\def\vv{\boldsymbol{v}}
\def\ww{\boldsymbol{w}}
\def\xx{\boldsymbol{x}}
\def\yy{\boldsymbol{y}}
\def\zz{\boldsymbol{z}}

\def\11{\mathbf{1}}

\def\00{\boldsymbol{0}}

\def\bbeta{\boldsymbol{\beta}}

\def\xxi{\boldsymbol{\xi}}

\def\frakX{\mathfrak{X}}
\def\frakD{\mathfrak{D}}

\def\transpose{\text{T}}

\def\spn{\operatorname{span}}

\DeclareMathAlphabet{\mathpzc}{OT1}{pzc}{m}{it}
\newcommand{\ce}{\mathfrak{e}}   
\newcommand{\cs}{\mathfrak{s}}   
\newcommand{\cf}{\mathfrak{f}}   

\DeclareMathOperator{\osc}{osc}
\DeclareMathOperator{\supp}{supp}
\DeclareMathOperator{\cond}{cond}
\DeclareMathOperator{\diag}{diag}

\title[Convergence of Adaptive Multilevel]{Convergence of Residual-driven Adaptive Kernel-based Multilevel Approximation}
\author{R\"udiger Kempf}
\address{Applied and Numerical Analysis, Department of Mathematics, University of Bayreuth, 95447 Bayreuth, Germany}
\email{ruediger.kempf@uni-bayreuth.de}
\date{}

\begin{document}

\begin{abstract}
Kernel-based multilevel methods combine compactly supported radial basis functions with a residual correction scheme and are known to converge linearly in the number of levels. In this paper we study an adaptive version of the multilevel method, in which the interpolation on each level only uses those sites at which the residual of the preceding level exceeds a given threshold. Although this method has been proposed some time ago and converges in all numerical tests, a convergence proof has been missing so far. We close this gap. We derive an a priori error bound, an a posteriori error bound valid for arbitrary point selections, and extend the analysis to penalized least-squares approximation and to target functions outside the native space of the kernel. For a variant with an additional check-loop per level, we show convergence with the rate of the classical multilevel method. We illustrate our findings with three numerical examples.
\end{abstract}



\maketitle

\section{Introduction}

Kernel-based methods offer a versatile tool for the reconstruction of an unknown function from scattered data. Their applications range from approximation \cite{Fasshauer2007,Wendland2004} to meshless numerical methods for partial differential equations \cite{Fornberg2015,Franke1998,Wendland1999} and machine learning \cite{Cucker2007,Schoelkopf2001}, most recently in application to operator learning \cite{Batlle2024,Kempf2026}. The great advantage they offer is the flexibility, stability, and deterministic error analysis \cite{Narcowich2005,NarcowichWardWendland2006,Schaback1995-2}.

If compactly supported radial basis functions are used, a trade-off between stability and accuracy has to be made: a fixed support radius leads to ill-conditioned systems if the number of data becomes large \cite{Schaback1995-2}, while a support radius coupled to the fill distance does not lead to convergence. This is known as the \emph{trade-off principle} \cite{Schaback1997}. The kernel-based multilevel method, first suggested in \cite{FloaterIske1996}, resolves this dilemma by computing a coarse approximation first and then successively correcting the residual on finer and finer sets of sites with kernels of smaller support, effectively representing multiple scales within the data. Although the method was established early and its convergence was demonstrated numerically, a proof of convergence for approximation on bounded domains is due to \cite{Wendland2010}, see also \cite{LeGiaSloanWendland2010} for the sphere and \cite{Wendland2017} for an overview. It is also used in collocation for partial differential equations \cite{Farrell2013} and high-dimensional problems \cite{BuettnerKempfWendland2026,Kempf2026_2,KempfWendland2023}. For globally supported kernels, the multilevel method can be realized efficiently in samplet coordinates \cite{AvesaniKempfMultererWendland2025}.

On many levels, however, the residual is already small in large parts of the domain. Hence, it is natural to use only those sites on which the residual is still non-negligible. 
Two such adaptive versions have been introduced in \cite{LeGiaWendland2014,Wendland2017}. The first one, discarding small coefficients after the computation of each local interpolant, goes back to \cite{LeGiaWendland2014} and its convergence is well-understood. The second one is a truly adaptive method: it evaluates the residual on the sites of the upcoming level and only uses those sites at which the absolute value of the residual exceeds a level-dependent threshold. This reduces the size of the linear systems and potentially the computational cost. Although all numerical tests indicate convergence, a proof has been missing so far. The difficulty is that the techniques of \cite{Wendland2010} require control of the residual on the full set of sites, whereas the adaptive method only enforces the residual to vanish on the selected subset, and we have no information about the residual on the remaining sites.

Recently, the idea of adaptive multilevel approximation has again gained interest and two related approaches have been proposed. The three approaches \cite{LotRieger2026}, \cite{AvesaniIskeMultererWenzel2026}, and \cite{Wendland2017} differ in what is selected and when. In \cite{LotRieger2026}, the sites are selected before the local approximation is computed, using the function values on the finest level, and the interpolant on the full set of sites is evaluated with the data at the removed sites set to zero. In \cite{AvesaniIskeMultererWenzel2026}, the centers of each level are chosen one by one by a greedy algorithm, which yields a bound on the number of centers. We follow the ideas of \cite{Wendland2017} described above.

The contribution of this paper is the first convergence analysis of the truly adaptive multilevel method of \cite{Wendland2017}. We derive an a posteriori error bound in $H^t(\Omega)$ and $L_\infty(\Omega)$, which holds for arbitrary selections of sites, in particular for the sets selected by the adaptive method itself. Besides the a priori term of the classical method, it only involves the residual at the dropped sites, measured in a discrete $\ell_2$-norm, and for the $L_\infty(\Omega)$ bound additionally its maximum on the dropped sites of the last level. These quantities can be evaluated at negligible cost for nested sites. Building on the same recursion, we prove an a priori error bound for the adaptive method, which shows linear convergence up to the thresholds. The main ingredients are the $\ell_2$-version of the sampling inequality, the observation that interpolation on an arbitrary subset of a quasi-uniform set of sites has uniformly exponentially localized Lagrange functions, and that the residual at the dropped sites is controlled by the smoothness of the residual of the previous level, which is measured on a coarser scale and therefore gains a factor of $\mu$, the refinement parameter of the sequence of sets of sites. This bound, however, requires the overlap parameter $\nu$ to be fixed before $\mu$ is chosen small, and its rate is of order $\sqrt{\mu}$ instead of the rate $\mu^s+\nu^{-s}$ of the classical method. Both restrictions disappear for a variant with an additional check-loop. We show convergence of this augmented adaptive method with the rate of the classical method, with constants independent of $\mu$ and $\nu$. This result also covers the coupled regime $\nu=\gamma/\mu$. We discuss why our technique does not yield this rate for the adaptive method itself, and conjecture that it holds. Finally, the analysis carries over to penalized least-squares approximation and to target functions outside the native space of the kernel. Numerical examples on two standard tests and one example function with distinct features show savings of up to a factor of eight in the number of centers compared with the classical method. Moreover, they show that the adaptive method attains the rate of the classical method, also in the coupled regime, which supports our conjecture.

The paper is organized as follows. In \cref{sec:multilevel} we recall the kernel-based multilevel method and, by slight modification in its proof, reduce the bound on the contraction factor that leads to the convergence of the method. In \cref{sec:adaptive} we introduce the two versions of the adaptive multilevel method, without and with check-loop, and discuss why the classical analysis fails. \Cref{sec:convergence,subsec:APosteriori,subsec:APriori} contain the error analysis for the version without check-loop, i.e., an a posteriori bound for arbitrary point selections and our main result, an a priori bound for the truly adaptive method. \cref{subsec:discussion} discusses the limitations of the analysis. In \cref{subsec:CheckLoopConvergence}, we prove convergence for the modified adaptive multilevel method with the check-loop and obtain the classical convergence rate. In \cref{sec:extensions} we extend these results to penalized least-squares approximation and to target functions which are rougher than the native space of the kernel. We give numerical examples in \cref{sec:numerics} and close with a summary in \cref{sec:summary}.

\section{The Kernel-based Multilevel Method}\label{sec:multilevel}
\subsection{The Setting}
Throughout this paper, $\Omega\subseteq\R^n$ denotes a bounded Lipschitz domain. For a finite set of sites $X := \{ \xx_1, \dots, \xx_N \} \subseteq\Omega$ we denote the \emph{fill distance} and \emph{separation radius} by
\begin{align*}
h_{X,\Omega}:=\sup_{\xx\in\Omega}\min_{\xx_j\in X}\|\xx-\xx_j\|_2,\quad \text{and} \quad q_X:=\frac12\min_{\xx_j\neq\xx_k}\|\xx_j-\xx_k\|_2.
\end{align*}
A sequence $ X_1, X_2, \dots \subseteq \Omega $ of sets of sites is called \emph{quasi-uniform} if there is a constant $c_{qu}\geq1$, independent of the level $ \ell $, such that $q_{X_{\ell}}\leq h_{X_{\ell},\Omega}\leq c_{qu}q_{X_{\ell}}$, $ \ell \in \N $.

We fix a continuous, compactly supported radial basis function $\Phi:\R^n\to\R$ with $\supp\Phi\subseteq\overline{B_1(\00)}$, which is a reproducing kernel of $H^s(\R^n)$, $s>n/2$, i.e., its Fourier transform satisfies
\begin{equation}\label{eq:fourier}
c_1\left(1+\|\xxi\|_2^2\right)^{-s}\leq\widehat{\Phi}(\xxi)\leq c_2\left(1+\|\xxi\|_2^2\right)^{-s}, \quad\xxi\in\R^n,
\end{equation}
with constants $0<c_1\leq c_2$, \cite{Wendland2004}. Every such $ \Phi $ gives rise to a \emph{native space} $ \cN_{\Phi}(\R^n) $, equipped with norm 
\begin{align*}
\| g \|_{\Phi}^2 := \int_{\R^n} \frac{|\widehat{g}(\xxi)|^2}{\widehat{\Phi}(\xxi)} \ d\xxi, \quad g \in \cN_{\Phi}(\R^n).
\end{align*}
For $\delta\in(0,1]$ the rescaled kernel $\Phi_\delta:=\delta^{-n}\Phi(\cdot/\delta)$ is again a reproducing kernel of $H^s(\R^n)$.
From \cite{Wendland2010} we have the norm equivalence
\begin{equation}\label{eq:normequiv}
c_1^{\frac{1}{2}}\|g\|_{\Phi_\delta}\leq\|g\|_{H^s(\R^n)}\leq c_2^{\frac{1}{2}}\delta^{-s}\|g\|_{\Phi_\delta},\quad g\in H^s(\R^n).
\end{equation}
Finally, $E:H^t(\Omega)\to H^t(\R^n)$ denotes the universal extension operator of Stein, which satisfies $Eu|_\Omega=u$ and $\|Eu\|_{H^t(\R^n)}\leq C_E\|u\|_{H^t(\Omega)}$ for all $0\leq t\leq s$, see \cite{Stein1970}. We note that the operator $ E $ is independent of $ t $, but the constant $ C_E $ depends on $ t $.

For the multilevel method we fix a sequence $X_1,X_2,\dots\subseteq\Omega$ of sets of sites with fill distances $h_{\ell} := h_{X_{\ell},\Omega}$ and separation radii $q_{\ell}:=q_{X_{\ell}}$. We assume that there is a \emph{refinement parameter} $\mu\in(0,1)$ such that $h_{\ell+1} \leq \mu h_{\ell}$, and we couple the scaling parameter to the fill distance by $\delta_{\ell} = \nu h_{\ell}$ with an \emph{overlap parameter} $\nu > 1$. We set $\Phi_{\ell} := \Phi_{\delta_{\ell}}$ and, for a set $\frakX \subseteq X_{\ell}$,
\begin{align*}
W_{\ell}(\frakX) := \spn \{ \Phi_{\ell}(\cdot-\xx):\xx\in \frakX \},
\end{align*}
with the convention $W_{\ell}(\emptyset)=\{0\}$. The interpolation operator $I_{\frakX,\Phi_{\ell}}:C(\Omega)\to W_{\ell}(\frakX)$ that maps $f$ to the unique interpolant of $f$ on $\frakX$ is given as 
\begin{align*}
I_{\frakX,\Phi_{\ell}}(f) = \sum_{\xx \in \frakX} \beta_{\xx} \Phi_{\ell}(\cdot - \xx),
\end{align*}
where the coefficient vector $ \bbeta \in \R^{\#\frakX} $ is the solution of the linear system $ M_{\frakX} \bbeta = f(\frakX) $, with $ M_{\frakX} = (\Phi_{\ell}(\xx - \zz))_{\xx,\zz \in \frakX} $.

It is well-known that $I_{\frakX,\Phi_{\ell}}f$ minimizes the native space norm among all interpolants in $ \cN_{\Phi_{\delta_{\ell}}} $, see \cite{Wendland2004}. For convenience we set $I_{\emptyset,\Phi_{\ell}}f := 0$. The classical multilevel method is a simple residual correction scheme: It starts with global interpolant $f_0=0$ and residual $e_0=f$, and, for $ \ell = 1, \dots, L $, it computes the local interpolants $s_{\ell}=I_{X_{\ell},\Phi_{\ell}}(e_{{\ell}-1}) \in W_{\ell}(X_{\ell}) $. Finally it updates the global interpolant $f_{\ell}=f_{{\ell}-1}+s_{\ell}$ and the residual $e_{\ell}=e_{{\ell}-1}-s_{\ell}$. The global interpolant $ f_L \in W_1(X_1) + \cdots + W_L(X_L) $ then interpolates $ f $ on $ X_L $.

For later reference, we fix the general assumptions.

\begin{assumption}\label{ass:standing}
We assume the following:
\begin{enumerate}
\item $\Omega \subseteq \R^n$ is a bounded Lipschitz domain.
\item $\Phi: \R^n \to \R$ is the reproducing kernel of $ H^s(\R^n)$, $ s > n/2$ and $ s \geq 1$, i.e. its Fourier transform satisfies \eqref{eq:fourier}, and $\supp \Phi \subseteq \overline{B_1(\00)}$.
\item The sets of sites $ X_{\ell} \subseteq \Omega$, ${\ell}\in\N$, with $ N_{\ell}$ elements, are quasi-uniform with a level-independent constant $c_{qu}$, and their fill distances satisfy $h_{{\ell}+1} \leq \mu h_{\ell}$ with $\mu \in (0,1)$ and $h_1$ sufficiently small.
\item The scaling parameters are given by $\delta_{\ell} = \nu h_{\ell} $ with $\nu>1$ and $\delta_1 \leq 1$.
\item The target function $ f $ is an element of $ H^s(\Omega) $.
\end{enumerate}
\end{assumption}

We do not assume that the sets of sites are nested. Nestedness is only relevant for the cost of evaluating the quantities $\rho_{\ell}$ and $\eta_{\ell}$ introduced in \cref{subsec:adaptiveSetting}. The assumption $s\geq1$ is only a restriction for $n=1$, and is always satisfied for Wendland's functions. The same holds for the condition $\sigma \geq 1$ in \cref{subsec:escaping}. The quasi-uniformity of the sets $ X_{\ell} $ means that 
\begin{align}\label{eq:Csharp}
C_{\#} := \sup_{\ell \in \N} h_{\ell}^n N_{\ell}
\end{align}
is finite.

The multilevel method is efficient. In the setting of \cref{ass:standing}, the number of nonzero entries per row and the condition number of the kernel matrices $ M_{X_{\ell}} $ are independent of the level. Hence, the number of steps of the conjugate gradient method without preconditioning is independent of the level as well, and the computational cost on each level is linear in the number of sites, see \cite[Corollary~10]{Wendland2017}.

\subsection{An Improved Convergence Result}

We now improve the convergence result of \cite[Theorem 10]{Wendland2017}, whose proof gives that there is a constant $ \widetilde{C} > 0 $ such that
\begin{align}\label{eq:classical}
\| E e_{\ell} \|_{\Phi_{\ell+1}} \leq \widetilde{C} (\nu^{-s} + \mu^s) \| Ee_{\ell - 1} \|_{\Phi_{\ell}}, \quad \ell \in \N.
\end{align}
In its proof, and throughout this paper, we need the following sampling inequalities, see \cite{ArcangeliLopezdeSilanesTorrens2007,WendlandRieger2005}. There are constants $h_0>0$, $C_t>0$, $0\leq t\leq s$, and $C_\infty>0$ such that for all $X\subseteq\Omega$ with $h:=h_{X,\Omega}\leq h_0$ and all $u\in H^s(\Omega)$
\begin{align}
\|u\|_{H^t(\Omega)}&\leq C_t\left(h^{s-t}\|u\|_{H^s(\Omega)}+h^{n/2-t}\|u\|_{\ell_2(X)}\right),\quad 0\leq t\leq s,\label{eq:sampling-Ht}\\
\|u\|_{L_\infty(\Omega)}&\leq C_\infty\left(h^{s-n/2}\|u\|_{H^s(\Omega)}+\|u\|_{\ell_\infty(X)}\right)\label{eq:sampling-Linf}
\end{align}
hold. We write $C_S:=C_0$ for the case $t=0$, i.e.,
\begin{align}\label{eq:sampling-l2}
\|u\|_{L_2(\Omega)}\leq C_S\left(h^{s}\|u\|_{H^s(\Omega)}+h^{n/2}\|u\|_{\ell_2(X)}\right),
\end{align}
and in particular
\begin{align}\label{eq:SamplingInequality}
\|u\|_{L_2(\Omega)}\leq C_S h^{s}\|u\|_{H^s(\Omega)}\quad\text{if } u|_X=0.
\end{align}
We also use the interpolation inequality, see, e.g., \cite{BrennerScott2008},
\begin{align}\label{eq:interpolation}
\|u\|_{H^1(\Omega)}\leq C_I\|u\|_{L_2(\Omega)}^{1-1/s}\|u\|_{H^s(\Omega)}^{1/s},\quad u\in H^s(\Omega).
\end{align}

Moreover, the residuals of the multilevel method satisfy, see \cite{Wendland2010},
\begin{align}\label{eq:NormEstimates}
\|e_{\ell}\|_{H^s(\Omega)}\leq c_2^{\frac{1}{2}}\delta_{\ell}^{-s}\|Ee_{\ell-1}\|_{\Phi_{\ell}}.
\end{align}

Our improvement is in that we introduce a variable scaling of the two summands in $ (\nu^{-s} + \mu^s)$ that we can optimize to reduce the contraction factor.

\begin{theorem}\label{thrm:KappaConvergence}
Let \cref{ass:standing} hold. Let $ \kappa > 0 $. Then the residuals $ e_{\ell} $ of the multilevel method satisfy
\begin{align}\label{eq:RecursionEstimateKappa}
\| Ee_{\ell} \|_{\Phi_{\ell+1}} \leq \alpha(\kappa) \| Ee_{\ell-1} \|_{\Phi_{\ell}},\quad \ell\in\N,
\end{align}
where
\begin{align}\label{eq:alphaKappa}
\alpha(\kappa) := C_K \max(C_S,1) \left((1 + \kappa^2)^{\frac{s}{2}} \nu^{-s} + (1 + \kappa^{-2})^{\frac{s}{2}} \mu^{s}\right),
\end{align}
with $ C_K := (c_2/c_1)^{1/2} C_E$, and $ C_S $ is the constant of \eqref{eq:SamplingInequality}.
\end{theorem}

\begin{proof}
With \eqref{eq:fourier} we have
\begin{align*}
\| Ee_{\ell} \|_{\Phi_{\ell+1}}^2 &\leq \frac{1}{c_1} \int_{\R^n} |\widehat{Ee_{\ell}}(\xxi)|^2 (1+ \delta_{\ell+1}^2 \| \xxi \|_2^2)^s \ d\xxi =: \frac{1}{c_1} (I_1 + I_2),
\end{align*}
where
\begin{align*}
I_1 = \int_{\| \xxi \|_2 \leq \frac{\kappa}{\delta_{\ell+1}}} |\widehat{Ee_{\ell}}(\xxi)|^2 (1+ \delta_{\ell+1}^2 \| \xxi \|_2^2)^s \ d\xxi
\end{align*}
and
\begin{align*}
I_2 = \int_{\| \xxi \|_2 \geq \frac{\kappa}{\delta_{\ell+1}}} |\widehat{Ee_{\ell}}(\xxi)|^2 (1+ \delta_{\ell+1}^2 \| \xxi \|_2^2)^s \ d\xxi.
\end{align*}
We emphasize that this is almost the same splitting as \cite[Proof Theorem 1]{Wendland2010}, with the important difference that we split along $ \| \xxi \|_2 = \frac{\kappa}{\delta_{\ell+1}}$, where \cite{Wendland2010} splits along $\| \xxi \|_2 = \frac{1}{\delta_{\ell+1}}$. The rest of the proof is virtually the same, however, for the convenience of the reader, we provide it.

To bound $ I_1 $, we note that we have $ (1 + \delta_{\ell+1}^2 \| \xxi \|_2^2)^s \leq (1 + \kappa^2)^s$, since $ \| \xxi \|_2^2 \delta_{\ell+1}^2 \leq \kappa^2 $. Using the properties of the extension operator, \eqref{eq:SamplingInequality}, and \eqref{eq:NormEstimates} yields
\begin{align*}
I_1 &\leq \left(1+\kappa^2\right)^s \|Ee_{\ell}\|_{L_2(\R^n)}^2 \\
&\leq \left(1+\kappa^2\right)^s C_E^2 \|e_{\ell}\|_{L_2(\Omega)}^2 \\
&\leq\left(1+\kappa^2\right)^s C_E^2 C_S^2 h_{\ell}^{2s} \|e_{\ell}\|_{H^s(\Omega)}^2\\
&\leq\left(1+\kappa^2\right)^sC_E^2C_S^2c_2\left(\frac{h_{\ell}}{\delta_{\ell}}\right)^{2s}\|Ee_{\ell-1}\|_{\Phi_{\ell}}^2 \\
&= \left(1+\kappa^2\right)^sC_E^2C_S^2c_2\,\nu^{-2s}\|Ee_{\ell-1}\|_{\Phi_{\ell}}^2.
\end{align*}

To bound $ I_2 $, we note that $ \| \xxi \|_2 \geq \kappa \delta_{\ell+1}^{-1} $ leads to $ 1 \leq \kappa^{-2} \delta_{\ell+1}^2 \| \xxi \|_2^2 $ and therefore
\begin{align*}
1+\delta_{\ell+1}^2\|\xxi\|_2^2\leq\left(1+\kappa^{-2}\right)\delta_{\ell+1}^2\|\xxi\|_2^2\leq\left(1+\kappa^{-2}\right)\delta_{\ell+1}^2\left(1+\|\xxi\|_2^2\right).
\end{align*}
Raising this to the power $s$, enlarging the domain of integration to $\R^n$, and using the boundedness of $E$ and \eqref{eq:NormEstimates} gives
\begin{align*}
I_2&\leq\left(1+\kappa^{-2}\right)^s\delta_{\ell+1}^{2s}\|Ee_{\ell}\|_{H^s(\R^n)}^2 \\
&\leq\left(1+\kappa^{-2}\right)^s\delta_{\ell+1}^{2s}C_E^2\|e_{\ell}\|_{H^s(\Omega)}^2\\
&\leq\left(1+\kappa^{-2}\right)^sC_E^2c_2\left(\frac{\delta_{\ell+1}}{\delta_{\ell}}\right)^{2s}\|Ee_{\ell-1}\|_{\Phi_{\ell}}^2 \\
&\leq \left( 1 + \kappa^{-2} \right)^s C_E^2 c_2\,\mu^{2s}\|Ee_{\ell-1}\|_{\Phi_{\ell}}^2.
\end{align*}
Adding both estimates, using $C_S \leq \max(C_S,1)$ and $\sqrt{a+b}\leq\sqrt a+\sqrt b$ yields \eqref{eq:RecursionEstimateKappa}.
\end{proof}

For $\kappa=1$, \eqref{eq:alphaKappa} is the contraction factor of \cite[Theorem 1]{Wendland2010} with $\widetilde{C} = 2^{s/2} C_K \max(C_S,1)$. This is the exact form of $ \widetilde{C} $ in \eqref{eq:classical}. 

The next goal is now to minimize $ \alpha(\kappa) $. To do this, we need the following elementary result.

\begin{lemma}\label{lem:minimization}
Let $a,b,q>0$ and $F(t):=b(1+t)^q+a(1+t^{-1})^q$ for $t>0$. Then $F$ attains its minimum at $t_*=(a/b)^{1/(q+1)}$, and
 \begin{align*}
\min_{t>0}F(t)=F(t_*)=\left(a^{1/(q+1)}+b^{1/(q+1)}\right)^{q+1}.
 \end{align*}
Moreover, $F(t_*)\leq F(1)=2^q(a+b)$ with equality if and only if $a=b$.
\end{lemma}

\begin{proof}
Since $ F $ is continuous and $F(t)\to\infty$ for $t\to0$ and $t\to\infty$, there is a minimum. Using $(1+t^{-1})^{q-1}=(1+t)^{q-1}t^{1-q}$ we obtain
\begin{align*}
F'(t)=qb(1+t)^{q-1}-qa\left(1+t^{-1}\right)^{q-1}t^{-2}=q(1+t)^{q-1}\left(b-at^{-(q+1)}\right).
\end{align*}
Hence, $t_*=(a/b)^{1/(q+1)}$ is the only critical point in $(0,\infty)$. Writing $F(t)=(1+t)^q(b+at^{-q})$ and using $at_*^{-q}=bt_*$ yields
\begin{align*}
F(t_*)=(1+t_*)^q\,b(1+t_*)=b\left(1+\left(\frac{a}{b}\right)^{\frac{1}{q+1}}\right)^{q+1}=\left(a^{\frac{1}{q+1}}+b^{\frac{1}{q+1}}\right)^{q+1}.
\end{align*}
The last statement follows since $t_*=1$ if and only if $a=b$, and $t_*$ is the unique minimizer.
\end{proof}

This allows us to find the smallest contraction parameter $ \alpha(\kappa) $.

\begin{corollary}\label{cor:kappa-opt}
With the notation and assumptions of \cref{thrm:KappaConvergence}, the function $\kappa\mapsto\alpha(\kappa)$ is minimized by $\kappa_*^2=(\mu\nu)^{2s/(s+2)}$. This yields $ \|Ee_{\ell}\|_{\Phi_{\ell+1}} \leq \alpha_* \|Ee_{\ell-1}\|_{\Phi_{\ell}}$, $\ell\in\N$, where
\begin{align*}
\alpha_*:=\alpha(\kappa_*) = C_K \max(C_S,1) \left(\mu^{\frac{2s}{s+2}} + \nu^{-\frac{2s}{s+2}} \right)^{\frac{s+2}{2}}.
\end{align*}
Moreover, $\alpha_* \leq \alpha(1) = 2^{s/2} C_K \max(C_S,1)(\mu^s+\nu^{-s})$ with equality if and only if $\mu\nu=1$.
\end{corollary}

\begin{proof}
We apply \cref{lem:minimization} with $t=\kappa^2$, $q=s/2$, $a=\mu^s$ and $b=\nu^{-s}$. This gives $t_*=(\mu^s\nu^s)^{2/(s+2)}$ and
\begin{align*}
\min_{\kappa>0}\left(\left(1+\kappa^2\right)^{\frac{s}{2}} \nu^{-s} + \left(1+\kappa^{-2}\right)^{\frac{s}{2}} \mu^s \right) = \left(\mu^{\frac{2s}{s+2}} + \nu^{-\frac{2s}{s+2}}\right)^{\frac{s+2}{2}}.
\end{align*}
The comparison with $\kappa=1$ is the last statement of \cref{lem:minimization}, since $a=b$ if and only if $\mu\nu=1$.
\end{proof}

The method in \cite{Wendland2010} considers a coupling of $\gamma:=\mu\nu$ with a fixed $ \gamma > 0 $. In this case, the bound of \cref{cor:kappa-opt} reads
\begin{align*}
\alpha_* = C_K \max(C_S,1) \left(1+\gamma^{-\frac{2s}{s+2}}\right)^{\frac{s+2}{2}} \,\mu^s.
\end{align*}
Compared with $ 2^{s/2} (1 + \gamma^{-s}) $ from \cite[(4.3.8)]{KempfThesis}, the factor $ (1 + \gamma^{-2s/(s+2)})^{(s+2)/2} $ is smaller unless $ \gamma = 1 $, in that case they are equal.

\section{The Adaptive Multilevel Method}\label{sec:adaptive}

\subsection{The Setting}\label{subsec:adaptiveSetting}

We now describe the truly adaptive multilevel method of \cite{Wendland2017}. In addition to the sets of sites $X_{\ell}$, which determine the parameters $h_{\ell}$, $\delta_{\ell}$ and $\Phi_{\ell}$, we construct adaptive sets of sites $\frakX_{\ell} \subseteq X_{\ell}$ which are actually used for the local interpolation. After the computation of the local interpolant on level $\ell$, we evaluate the residual on $X_{\ell+1}$ and keep only those sites at which its absolute value exceeds a given threshold $\varepsilon_{\ell} \geq 0$. We note that this comes at virtually no additional cost, since the standard multilevel method also has to evaluate $ e_{\ell} $ on $ X_{{\ell}+1} $ to obtain the data for the interpolation on level $ \ell + 1 $. To distinguish the quantities of the adaptive method from those of the classical multilevel method, we denote them by fraktur letters: $\cs_{\ell}$ is the adaptive local interpolant, $\ce_{\ell}$ the adaptive residual and $\cf_{\ell}$ the adaptive global approximation. Note that $\ce_{\ell}$ and $e_{\ell}$ differ in general, since $\cs_{\ell}$ only interpolates on $\frakX_{\ell}$. The method is summarized in \cref{alg:adaptive}.

\begin{algorithm}[ht]
\caption{Adaptive multilevel approximation}\label{alg:adaptive}
\KwData{number of levels $L$, function $f$, thresholds $\varepsilon_{\ell} \geq 0$, sets of sites $X_1,\ldots,X_L$}
\KwResult{approximation $\cf_L \in W_1(\frakX_1) + \cdots + W_L(\frakX_L)$}
Set $\cf_0=0$, $\ce_0=f$, $\frakX_1=X_1$\;
\For{$\ell = 1,2,\ldots,L$}{
    Compute $\cs_{\ell} = I_{\frakX_{\ell},\Phi_{\ell}}(\ce_{{\ell}-1})$\;
  Set $\cf_{\ell} = \cf_{{\ell}-1} + \cs_{\ell}$ and $\ce_{\ell} = \ce_{{\ell} - 1} -\cs_{\ell}$\;
 \If{$\ell < L$}{
    Set $\frakX_{{\ell}+1} = \{\xx \in X_{{\ell}+1} \;: \; |\ce_{\ell}(\xx)| > \varepsilon_{\ell} \}$\;
  }
}
\end{algorithm}

If $\frakX_{\ell} = \emptyset$, the residual is bounded by $\varepsilon_{{\ell} - 1}$ on all of $X_{\ell}$, no linear system has to be solved, and level ${\ell}$ is skipped, i.e., $\ce_{\ell} = \ce_{{\ell} - 1}$. The algorithm then continues with the next level, on which the residual is checked on the finer set $X_{{\ell}+1}$. All results below include this case.

For the analysis we denote the dropped sites on level ${\ell}$ by $\frakD_{\ell} := X_{\ell} \setminus \frakX_{\ell}$. Clearly, $\frakD_1=\emptyset$ and, by construction,
\begin{align}\label{eq:G}
|\ce_{{\ell}-1}(\xx)| \leq \varepsilon_{{\ell}-1},\quad\xx\in\frakD_{\ell},\quad \ell \geq 2.
\end{align}

Since the local interpolant on level $\ell$ only interpolates on $\frakX_{\ell}$, the residual $\ce_{\ell}$ vanishes on $ \frakX_{\ell} $, but on $\frakD_{\ell}$ we have no direct information, since
\begin{align*}
\ce_{\ell}(\xx) = \ce_{{\ell}-1}(\xx) - \cs_{\ell}(\xx), \quad \xx \in \frakD_{\ell},
\end{align*}
and $\cs_{\ell}(\xx)$ is not controlled by the threshold. This is the reason why \eqref{eq:SamplingInequality} and hence the proof of \cite[Theorem 1]{Wendland2010} cannot be applied directly. One might try to use the sampling inequality
\begin{align*}
\|u\|_{L_2(\Omega)}\leq C\left(h_{X,\Omega}^s\|u\|_{H^s(\Omega)}+\|u\|_{\ell_\infty(X)}\right)
\end{align*}
instead, which requires bounding $\|\ce_{\ell}\|_{\ell_\infty(\frakD_{\ell})}$. This, however, does not lead to a convergence result. The bound $\|g\|_{L_\infty(\R^n)}\leq C\delta_{\ell}^{-n/2}\|g\|_{\Phi_{\ell}}$ is sharp, so that an $\ell_\infty$-term cannot be compared with the native space norm without losing a factor $\delta_{\ell}^{-n/2}$. Instead, we measure the residual on the dropped sites with the discrete $\ell_2$-norm and define, for $\ell \in \N$,
\begin{align}\label{eq:rho}
\rho_{\ell} := h_{\ell}^{\frac{n}{2}} \left( \sum_{\xx \in \frakD_{\ell}} |\ce_{\ell}(\xx)|^2\right)^{\frac{1}{2}}.
\end{align}
Since $\ce_{\ell}$ vanishes on $\frakX_{\ell}$, we also have $\rho_{\ell} = h_{\ell}^{n/2} \|\ce_{\ell}\|_{\ell_2(X_{\ell})}$. Note that $\rho_1 = 0$ and that $\rho_{\ell}$ can be computed with $\cO(\#X_{\ell})$ operations.

We also set $\eta_\ell := \|\ce_\ell\|_{\ell_\infty(\frakD_\ell)}$, with $\eta_\ell := 0$ if $\frakD_\ell=\emptyset$.

In \cite{Wendland2017} a variation of \cref{alg:adaptive} was proposed: By including an additional check-loop on each level, it is possible to guarantee that, once the local interpolant has been computed, the residual stays below a second threshold $\tau_{\ell}$ on all sites of the current level, including the dropped ones. To this end, all dropped sites at which the residual exceeds $\tau_{\ell}$ are added to the selected set, and the local interpolant is recomputed until no such site remains. The method is summarized in \cref{alg:checkloop}.

\begin{algorithm}[ht]
\caption{Adaptive multilevel approximation with check-loop}\label{alg:checkloop}
\KwData{number of levels $L$, function $f$, thresholds $\varepsilon_{\ell} $ and $\tau_{\ell} \geq 0$, sets of sites $X_1,\ldots,X_L$}
\KwResult{approximation $\cf_L \in W_1(\frakX_1) + \cdots + W_L(\frakX_L)$}
Set $\cf_0=0$, $\ce_0=f$, $\frakX_1=X_1$\;
\For{$\ell = 1,2,\ldots,L$}{
  Compute $\cs_{\ell} = I_{\frakX_{\ell},\Phi_{\ell}}(\ce_{{\ell}-1})$ and $\ce_{\ell} = \ce_{{\ell} - 1} -\cs_{\ell}$\;
  \While{$|\ce_{\ell}(\xx)| > \tau_{\ell}$ for some $\xx \in \frakD_{\ell} = X_{\ell} \setminus \frakX_{\ell}$}{
    Add all $\xx \in \frakD_{\ell}$ with $|\ce_{\ell}(\xx)| > \tau_{\ell}$ to $\frakX_{\ell}$\;
    Recompute $\cs_{\ell} = I_{\frakX_{\ell},\Phi_{\ell}}(\ce_{{\ell}-1})$ and $\ce_{\ell} = \ce_{{\ell} - 1} -\cs_{\ell}$\;
  }
  Set $\cf_{\ell} = \cf_{{\ell}-1} + \cs_{\ell}$\;
  \If{$\ell < L$}{
    Set $\frakX_{{\ell}+1} = \{\xx \in X_{{\ell}+1} \;: \; |\ce_{\ell}(\xx)| > \varepsilon_{\ell} \}$\;
  }
}
\end{algorithm}

It is clear that each loop of the algorithm terminates after at most $ \# X_{\ell} $ steps.

\subsection{Auxiliary Results}\label{sec:convergence}

For later reference, we collect known statements, adapted for $ \ce_{\ell} $. 

\begin{lemma}\label{lem:stability}
Let \cref{ass:standing} hold. For every $\ell \in \N$ and every choice of $\frakX_{\ell} \subseteq X_{\ell}$ the estimates
\begin{align}\label{eq:stab-a}
\|\ce_{\ell}\|_{H^s(\Omega)} &\leq c_2^{\frac{1}{2}} \delta_{\ell}^{-s} \| E\ce_{{\ell} - 1}\|_{\Phi_{\ell}},
\end{align}
and 
\begin{align}\label{eq:stab-b}
\|\ce_{\ell}\|_{L_2(\Omega)} &\leq C_S \left(c_2^{\frac{1}{2}} \nu^{-s} \|E\ce_{{\ell}-1}\|_{\Phi_{\ell}} + \rho_{\ell}\right)
\end{align}
hold.
\end{lemma}

\begin{proof}
Since $\cs_\ell$ interpolates $E\ce_{\ell-1}$ on $\frakX_\ell$, the function $E\ce_{\ell-1}-\cs_\ell$
extends $\ce_\ell$, and $\cs_\ell$ is the $\Phi_\ell$-orthogonal projection of $E\ce_{\ell-1}$
onto $W_\ell(\frakX_\ell)$, also for $\frakX_\ell=\emptyset$. With \eqref{eq:normequiv} and $\delta_\ell \leq 1$ this gives
\begin{align*}
\|\ce_\ell\|_{H^s(\Omega)} \leq \|E\ce_{\ell-1}-\cs_\ell\|_{H^s(\R^n)}
\leq c_2^{1/2}\delta_\ell^{-s}\|E\ce_{\ell-1}-\cs_\ell\|_{\Phi_\ell}
\leq c_2^{1/2}\delta_\ell^{-s}\|E\ce_{\ell-1}\|_{\Phi_\ell}.
\end{align*}

 To prove \eqref{eq:stab-b}, we apply \eqref{eq:sampling-l2} to $\ce_{\ell}$ with the full set of sites $X_{\ell}$. Since $\ce_{\ell}$ vanishes on $\frakX_{\ell}$, this yields together with \eqref{eq:stab-a}
\begin{align*}
\|\ce_{\ell}\|_{L_2(\Omega)} \leq C_S \left(h_{\ell}^s\|\ce_{\ell}\|_{H^s(\Omega)} + h_{\ell}^{\frac{n}{2}} \| \ce_{\ell} \|_{\ell_2(X_{\ell})} \right) 
\leq C_S\left(c_2^{\frac{1}{2}} \left(\frac{h_{\ell}}{\delta_{\ell}}\right)^s \|E\ce_{{\ell}-1}\|_{\Phi_{\ell}} + \rho_{\ell} \right),
\end{align*}
which is \eqref{eq:stab-b}, since $h_{\ell}/\delta_{\ell}=\nu^{-1}$.
\end{proof}

With \cref{lem:stability} we can now follow the proof of \cite[Theorem 1]{Wendland2010} and obtain a perturbed recursion estimate.

\begin{lemma}\label{lem:recursion}
Let \cref{ass:standing} hold. For every $\ell \in \N$ and every choice of $\frakX_{\ell} \subseteq X_{\ell}$ the estimate
\begin{align}\label{eq:recursion}
\|E\ce_{\ell}\|_{\Phi_{\ell+1}} \leq \alpha \,\|E\ce_{\ell-1}\|_{\Phi_{\ell}} + C_P \,\rho_{\ell}
\end{align}
holds with $\alpha:=\widetilde{C}(\mu^s+\nu^{-s})$, where $\widetilde{C}$ is the constant of \eqref{eq:classical}, and
\begin{align*}
C_P:=\left(\frac{2^s}{c_1}\right)^{1/2}C_EC_S.
\end{align*}
\end{lemma}

\begin{proof}
In the proof of \cref{thrm:KappaConvergence}, the residual enters only through $\|\ce_\ell\|_{L_2(\Omega)}$
in $I_1$ and $\|\ce_\ell\|_{H^s(\Omega)}$ in $I_2$. Bounding these with \eqref{eq:stab-b} and \eqref{eq:stab-a}
instead of \eqref{eq:SamplingInequality} and \eqref{eq:NormEstimates} gives
\begin{align*}
\|E\ce_\ell\|_{\Phi_{\ell+1}} \leq \alpha(\kappa) \|E\ce_{\ell-1}\|_{\Phi_\ell} + C_P(\kappa)\rho_\ell,
\quad \text{where} \quad C_P(\kappa) := \left( \frac{(1+\kappa^2)^s}{c_1} \right)^{1/2} C_E C_S,
\end{align*}
for every $\kappa>0$. For $\kappa=1$ this is \eqref{eq:recursion}.
\end{proof}

\begin{remark}\label{rem:kappa-adaptive}
 All results below hold for every $\kappa>0$ with $\alpha$ and $C_P$ replaced by $\alpha(\kappa)$ and $C_P(\kappa)$, in particular with the optimal contraction parameter $\alpha_*=\alpha(\kappa_*)$ from \cref{cor:kappa-opt}. Note that $C_P(\kappa_*)\leq C_P$ if and only if $\mu\nu\leq1$. For the sake of readability, we state all results with $\alpha$ and $C_P$.
\end{remark}

The next result shows that interpolation on an arbitrary subset of $X_{\ell}$ is uniformly localized. For the full set $X_{\ell}$, the exponential decay of the Lagrange functions of a stationary rescaled kernel and the boundedness of the Lebesgue constant are shown in \cite[Theorem~2.3, Corollary~2.5]{DeMarchiWendland2020}. The proof is based on the exponential decay of the entries of the inverse of a banded matrix \cite[Theorem~3.3]{BaxterSivakumarWard1994}, see also \cite{Demko1984,NarcowichWard1991}, and it only uses the separation of the sites, the number of sites in a ball of radius $\delta_{\ell}$, and upper bounds on the norm of the inverse and on the condition number of the kernel matrix. All of these hold uniformly for subsets: the separation radius of $\frakX \subseteq X_{\ell} $ is bounded from below by $q_{\ell}$, and since the kernel matrix of $\frakX$ is a principal submatrix of the one of $X_{\ell}$, the interlacing of eigenvalues shows that deleting sites cannot increase its condition number. The fill distance of $\frakX$, which may be much larger than the one of $ X_{\ell} $, e.g., if $\frakX$ has large holes, is not needed. For a nonempty set $\frakX \subseteq X_{\ell}$ we denote by $L_{\xx}^{\frakX}\in W_{\ell}(\frakX)$, $\xx\in\frakX$, the Lagrange functions, defined by $L_{\xx}^{\frakX}(\zz)=\delta_{\xx,\zz}$ for $\xx,\zz\in\frakX$.

\begin{lemma}\label{lem:localization}
Let \cref{ass:standing} hold. There are constants $C_L>0$ and $\eta>0$, depending only on $\Phi$, $n$, $s$, $\nu$ and $c_{qu}$, such that for every $\ell\in\N$ and every nonempty set $\frakX\subseteq X_{\ell}$ the estimate
\begin{align*}
|L_{\xx}^{\frakX}(\yy)|\leq C_L\exp\left(-\eta\frac{\|\xx-\yy\|_2}{q_{\ell}}\right)\leq C_L\,\theta^{\|\xx-\yy\|_2/\delta_{\ell}}, \quad \xx\in\frakX,\ \yy\in\R^n,
\end{align*}
holds, where $\theta:=e^{-\eta\nu}\in(0,1)$.
\end{lemma}

\begin{proof}
Let $M_{\frakX}=(\Phi_{\ell}(\xx-\zz))_{\xx,\zz\in\frakX}$ denote the kernel matrix on $\frakX$. It is a principal submatrix of the kernel matrix $M_{\ell}$ on $X_{\ell}$. For $ \vv \in \R^{\#\frakX} $ let $\widetilde{\vv} \in \R^{\#X_{\ell}} $ denote its extension by zero. Then $\vv^\transpose M_{\frakX} \vv = \widetilde{\vv}^\transpose M_{\ell} \widetilde{\vv}$ and $ \| \vv \|_2 = \| \widetilde{\vv} \|_2$, and the Courant-Fischer theorem yields
\begin{align*}
\lambda_{\min}(M_{\frakX})=\min_{\vv\neq\00}\frac{\vv^\transpose M_{\frakX}\vv}{\vv^\transpose\vv}\geq\lambda_{\min}(M_{\ell}), \quad \text{and} \quad \lambda_{\max}(M_{\frakX})\leq\lambda_{\max}(M_{\ell}).
\end{align*}
This is a special case of Cauchy's interlacing theorem, see, e.g., \cite[Corollary~III.1.5]{Bhatia1997}. This yields
\begin{equation}\label{eq:cond}
\begin{aligned}
\|M_{\frakX}^{-1}\|_2&\leq\|M_{\ell}^{-1}\|_2\leq C_A(\nu c_{qu})^{2s-n}\delta_{\ell}^n, \quad \text{and } \\
\cond_2(M_{\frakX})&\leq \cond_2(M_{\ell}) \leq C_A(\nu c_{qu})^{2s-n}(1+\nu c_{qu})^n\Phi(\00) =: \kappa_\nu,
\end{aligned}
\end{equation}
with a constant $C_A$ depending only on $\Phi$, $n$ and $s$. 

Since $ \kappa_{\nu} $ is an upper bound on the condition number, we may replace it with $ \max(\kappa_{\nu},4) $ and use $\kappa_\nu \geq 4$. This ensures that $ \eta $ below is well defined, also if $ \#\frakX = 1 $.

The sites of $\frakX$ are $2q_{\ell}$-separated. Hence, we can relabel $M_{\frakX}$ as in \cite[Section~2]{DeMarchiWendland2020} as an $R$-banded multivariate matrix with $ R := \max(\nu c_{qu},4)\sqrt n$. \Cref{eq:cond} and \cite[Theorem~3.3]{BaxterSivakumarWard1994} yield as in the proof of \cite[Theorem~2.3]{DeMarchiWendland2020}
\begin{align*}
|(M_{\frakX}^{-1})_{\xx,\zz}| \leq 2\|M_{\frakX}^{-1}\|_2\,e^{2\eta\sqrt n}\exp\left(-\eta\frac{\|\xx-\zz\|_2}{q_{\ell}}\right), \quad\xx,\zz\in\frakX,
\end{align*}
with $ \eta:=\frac{1}{2R}\log\frac{\sqrt{\kappa_\nu}+1}{\sqrt{\kappa_\nu}-1} $. Now let $\yy\in\R^n$. The Lagrange functions have the representation
\begin{align*}
L_{\xx}^{\frakX}(\yy) = \sum_{\zz\in\frakX}(M_{\frakX}^{-1})_{\xx,\zz}\Phi_{\ell}(\yy-\zz),
\end{align*}
in which only the sites $\zz$ with $ \| \yy - \zz \|_2 < \delta_{\ell}$ contribute because of the compact support of $ \Phi_{\ell} $. By \cite[(9)]{DeMarchiWendland2020} there are at most $(1+\nu c_{qu})^n$ of them, and for each of them we have $|\Phi_{\ell}(\yy-\zz)|\leq\delta_{\ell}^{-n}\Phi(\00)$ and $\|\xx-\zz\|_2\geq\|\xx-\yy\|_2-\delta_{\ell}$. Putting everything together and using $\delta_{\ell}/q_{\ell}\leq\nu c_{qu}$ gives the first estimate with
\begin{align*}
C_L:=2C_A(\nu c_{qu})^{2s-n}(1+\nu c_{qu})^n\Phi(\00)\,e^{\eta(2\sqrt n+\nu c_{qu})}.
\end{align*}
The second estimate follows from $ q_{\ell} \leq h_{\ell}=\delta_{\ell}/\nu$.
\end{proof}

For $m\in\N_0$ we define the quantities
\begin{align*}
b_m:=C_L\theta^m\left(1+(m+1)\nu c_{qu}\right)^n,\quad \Lambda := \sum_{m=0}^\infty b_m,\quad \text{and} \quad B_s:=\sum_{m=0}^\infty b_m(m+1)^{2s},
\end{align*}
and, for $\yy\in X_{\ell}$, the annuli $A_m(\yy) := \{ \xx \in \R^n \, : \, m \delta_{\ell} \leq \| \xx - \yy \|_2 < (m+1) \delta_{\ell} \}$. A comparison of volumes shows $\#(X_{\ell} \cap A_m(\yy)) \leq \#(X_{\ell} \cap B_{(m+1) \delta_{\ell}}(\yy)) \leq (1+(m+1)\nu c_{qu})^n $, which was already used in the proof above. On $A_m(\yy)$ \cref{lem:localization} gives $|L_{\xx}^{\frakX}(\yy)|\leq C_L\theta^m$. Hence,
\begin{equation}\label{eq:annulus}
\sum_{\xx\in\frakX\cap A_m(\yy)}|L_{\xx}^{\frakX}(\yy)| \leq b_m, \quad\yy\in X_{\ell},\ m\in\N_0.
\end{equation}
In particular, $\sum_{\xx\in\frakX}|L_{\xx}^{\frakX}(\yy)|\leq\Lambda$ for all $\yy\in X_{\ell}$, i.e., the Lebesgue constant of interpolation on an arbitrary subset of $X_{\ell}$ is bounded independently of the level and of the subset. This is the analogue of \cite[Corollary~2.5]{DeMarchiWendland2020}.

\begin{remark}\label{rem:nu-localization}
The proof of \cref{lem:localization} shows how the constants depend on $\nu$. We are particularly interested in their behavior as $ \nu $ grows. To simplify, assume $ \nu $ is large enough such that $\nu c_{qu}\geq4$. This means that $ R = \nu c_{qu} \sqrt n$, and that $\kappa_\nu\leq C(\nu c_{qu})^{2s}$. Since $2/t\leq\log((t+1)/(t-1))\leq2/(t-1)$ for $t>1$ and $\kappa_\nu\geq4$, we have
\begin{align*}
\frac{c}{(\nu c_{qu})^{s+1}}\leq\frac{1}{R\sqrt{\kappa_\nu}}\leq\eta\leq\frac{1}{R(\sqrt{\kappa_\nu}-1)}\leq\frac{1}{R}.
\end{align*}
In particular, $\eta(2\sqrt n+\nu c_{qu})\leq2$, and $C_L\leq C(\nu c_{qu})^{2s}$. This means, since $\eta\nu\leq1$, also $1-\theta\geq\eta\nu/2\geq c\,\nu^{-s}$. The constants $c,C>0$ only depend on $\Phi$, $n$, $s$ and $c_{qu}$. Hence, the decay rate deteriorates as $\nu$ grows, since $\theta \to 1$. Using $\sum_{m\geq0}(m+2)^p\theta^m\leq C_p(1-\theta)^{-(p+1)}$, we obtain
\begin{align*}
\Lambda\leq C \nu^{2s+n+s(n+1)},\quad \text{and} \quad B_s \leq C\nu^{2s+n+s(n+2s+1)}.
\end{align*}
These bounds are far from being sharp. However, a polynomial growth of $\Lambda$ in $\nu$ is to be expected for subsets $\frakX$ with holes which are larger than $\delta_{\ell}$, since the Lagrange functions then extrapolate over a distance of up to $\delta_{\ell}=\nu h_{\ell}$ beyond the boundary of the hole.
\end{remark}

The next lemma is a discrete version of the fact that a function in $H^s(\R^n)$ oscillates only little on small balls.

\begin{lemma}\label{lem:oscillation}
Let \cref{ass:standing} hold. Let $ g \in H^s(\R^n)$ and $ r \geq h_{\ell}$. For $ \yy \in \R^n$ let
\begin{align*}
\omega_r(\yy):=\sup_{\|\zz-\yy\|_2\leq r}|g(\zz)-g(\yy)|.
\end{align*}
Then there is a constant $C_O>0$, depending only on $n$, $s$ and $c_{qu}$, such that
\begin{align*}
h_{\ell}^n \sum_{\yy\in X_{\ell}} \omega_r(\yy)^2 \leq C_O\left( r^2 |g|_{H^1(\R^n)}^2+r^{2s} |g|_{H^s(\R^n)}^2 \right).
\end{align*}
\end{lemma}

\begin{proof}
We first consider the unit cube $Q=(0,1)^n$ and $u\in H^s(Q)$. Let $\bar u$ denote the mean value of $u$ on $Q$. We have
\begin{align*}
\osc_Q u := \sup_{\zz,\yy\in Q} | u(\zz) - u(\yy)| \leq \sup_{\zz\in Q} |u(\zz)-\bar u| + \sup_{\yy\in Q}|\bar u-u(\yy)| = 2\|u-\bar u\|_{L_\infty(Q)}.
\end{align*}
Since $s>n/2$, the Sobolev embedding theorem \cite[Theorem~4.12]{AdamsFournier2003} gives $\|u-\bar u\|_{L_\infty(Q)}\leq C\|u-\bar u\|_{H^s(Q)}$. On the Lipschitz domain $Q$ the full norm is equivalent to the sum of the $L_2$-norm and the highest seminorm, i.e., $\|v\|_{H^s(Q)}\leq C(\|v\|_{L_2(Q)}+|v|_{H^s(Q)})$, see, e.g., \cite[Chapter~5]{AdamsFournier2003}. We apply this to $v=u-\bar u$ and note that $|u-\bar u|_{H^s(Q)} = |u|_{H^s(Q)}$, since the seminorm vanishes on constants. Finally, Poincar\'e's inequality \cite[Section~5.3]{BrennerScott2008} gives $\|u-\bar u\|_{L_2(Q)}\leq C|u|_{H^1(Q)}$. Altogether, we obtain
\begin{align*}
\osc_Qu\leq C\left(|u|_{H^1(Q)}+|u|_{H^s(Q)}\right).
\end{align*}
For a cube $Q_a$ of side length $a>0$ we apply this estimate to $u(a\,\cdot)$ and use $|u(a\,\cdot)|_{H^t(Q)}=a^{t-n/2}|u|_{H^t(Q_a)}$, which yields
\begin{align*}
\osc_{Q_a} u \leq C\left(a^{1-n/2}|u|_{H^1(Q_a)} + a^{s-n/2}|u|_{H^s(Q_a)}\right).
\end{align*}
Next, we tile $ \R^n $ with cubes $Q$ of side length $r$ and denote by $Q^*$ the concentric cube with side length $3r$. For $\yy\in Q$ we have $\overline{B_r(\yy)}\subseteq Q^*$ and hence $\omega_r(\yy)\leq\osc_{Q^*}g$. Since the sites of $X_{\ell}$ are separated, a comparison of volumes yields $\#(X_{\ell}\cap Q)\leq C_n(r/q_{\ell})^n\leq C_n(c_{qu}r/h_{\ell})^n$, where we used $r\geq h_{\ell}\geq q_{\ell}$. This gives
 \begin{align*}
h_{\ell}^n\sum_{\yy\in X_{\ell}\cap Q}\omega_r(\yy)^2\leq C_nc_{qu}^nr^n\left(\osc_{Q^*}g\right)^2\leq C\left(r^2|g|_{H^1(Q^*)}^2+r^{2s}|g|_{H^s(Q^*)}^2\right).
 \end{align*}
Finally, we sum over all cubes $Q$ and note that every point lies in at most $3^n$ of the cubes $Q^*$.
\end{proof}

\subsection{An A Posteriori Error Bound}\label{subsec:APosteriori}

Iterating the recursion \eqref{eq:recursion} and using a sampling inequality on the last
level yields an error bound for arbitrary point selections.

\begin{theorem}\label{thm:aposteriori}
Let \cref{ass:standing} hold. Let $\frakX_\ell\subseteq X_\ell$, $\ell\in\N$, be arbitrary
and let $\cf_L$ be the approximation computed by \cref{alg:adaptive} with the sets
$\frakX_\ell$. Then there exists a constant $C>0$, independent of $f$, $L$, $\mu$, $\nu$ and the
choice of the sets $\frakX_\ell$, such that for $0 \leq t \leq s$
\begin{align}\label{eq:aposteriori-Ht}
\| f - \cf_L \|_{H^t(\Omega)} \leq C h_L^{-t} \left( \alpha^L \|f\|_{H^s(\Omega)}
  + \sum_{k=1}^{L} \alpha^{L-k} \rho_k \right)
\end{align}
and
\begin{equation}\label{eq:aposteriori-Linf}
\|f-\cf_L\|_{L_\infty(\Omega)} \leq C \left( h_L^{-n/2} \left( \alpha^L \|f\|_{H^s(\Omega)}
  + \sum_{k=1}^{L-1} \alpha^{L-k} \rho_k \right) + \eta_L \right)
\end{equation}
hold for all $L\in\N$.
\end{theorem}

\begin{proof}
Applying \eqref{eq:recursion} $L-1$ times and using $\ce_0=f$ gives
\begin{align*}
\|E\ce_{L-1}\|_{\Phi_L} \leq \alpha^{L-1}\|Ef\|_{\Phi_1}
  + C_P\sum_{k=1}^{L-1}\alpha^{L-1-k}\rho_k,
\end{align*}
where the sum is empty for $L=1$. By \eqref{eq:normequiv} with $\delta_1 \leq 1$ and the boundedness of $E$, we
have $\|Ef\|_{\Phi_1} \leq c_1^{-1/2}C_E\|f\|_{H^s(\Omega)}$. Since
$\alpha \geq \widetilde{C}\nu^{-s}$, we have $\nu^{-s}\alpha^{L-1-k} \leq \widetilde{C}^{-1}\alpha^{L-k}$ for
$0 \leq k \leq L-1$, and therefore
\begin{align}\label{eq:aposteriori-aux}
\nu^{-s}\|E\ce_{L-1}\|_{\Phi_L} \leq C' \left( \alpha^L \|f\|_{H^s(\Omega)}
  + \sum_{k=1}^{L-1} \alpha^{L-k} \rho_k \right),
\end{align}
with $ C' := \widetilde{C}^{-1}\max\{c_1^{-1/2}C_E,\ C_P\} $.

We apply \eqref{eq:sampling-Ht} and \eqref{eq:sampling-Linf} to $u=\ce_L = f-\cf_L$.
By \eqref{eq:stab-a} and $h_L/\delta_L=\nu^{-1}$, we have
\begin{align*}
h_L^{s-t}\|\ce_L\|_{H^s(\Omega)} \leq c_2^{1/2}h_L^{-t}\nu^{-s}\|E\ce_{L-1}\|_{\Phi_L},
\quad 0 \leq t \leq s.
\end{align*}
Since $\ce_L$ vanishes on $\frakX_L$, we have $h_L^{n/2-t}\|\ce_L\|_{\ell_2(X_L)} = h_L^{-t}\rho_L$
and $\|\ce_L\|_{\ell_\infty(X_L)} = \eta_L$. Together with \eqref{eq:aposteriori-aux}, this
proves \eqref{eq:aposteriori-Ht} and \eqref{eq:aposteriori-Linf} with
$C := \max\{C_t,C_\infty\}\max\{c_2^{1/2}C',1\}$, which is independent of $f$, $L$, $\mu$,
$\nu$ and the sets $\frakX_\ell$.
\end{proof}

The quantities $\rho_k$ and $\eta_k$ can be computed a posteriori with $\mathcal{O}(\#X_k)$ operations.
If the sets of sites are nested, they come at no additional cost, since \cref{alg:adaptive} evaluates $\ce_k$ on $X_{k+1}\supseteq X_k$ anyway.

\subsection{An A Priori Error Bound}\label{subsec:APriori}

It remains to bound $\rho_{\ell}$ for \cref{alg:adaptive} itself. We start by bounding the residual on the dropped sites by the threshold and by the oscillation of the residual of the preceding level.

\begin{proposition}\label{prop:rho}
Let \cref{ass:standing} hold and let the sets $\frakX_{\ell}$ be produced by \cref{alg:adaptive}. Then the estimate
 \begin{align*}
\rho_{\ell}\leq K_1\varepsilon_{\ell-1}+K_2\left(\delta_{\ell}|E\ce_{\ell-1}|_{H^1(\R^n)}+\delta_{\ell}^s|E\ce_{\ell-1}|_{H^s(\R^n)}\right)
 \end{align*}
holds for all $\ell\geq2$ with $K_1:=(2C_\#)^{1/2}(1+\Lambda)$ and $K_2:=(2\Lambda C_OB_s)^{1/2}$.
\end{proposition}

\begin{proof}
We set $g:=E\ce_{\ell-1}$ and note that $g=\ce_{\ell-1}$ on $X_{\ell}\subseteq\Omega$. Let $\yy\in\frakD_{\ell}$. If $\frakX_{\ell}=\emptyset$, we have $\ce_{\ell}(\yy)=\ce_{\ell-1}(\yy)$ and hence $|\ce_{\ell}(\yy)|\leq\varepsilon_{\ell-1}$ by \eqref{eq:G}. Otherwise, the Lagrange representation of $\cs_{\ell}$ gives
 \begin{align*}
\ce_{\ell}(\yy) = \ce_{\ell-1}(\yy) - \cs_{\ell}(\yy) = \ce_{\ell-1}(\yy)-\sum_{\xx\in\frakX_{\ell}}L_{\xx}^{\frakX_{\ell}}(\yy)\,\ce_{\ell-1}(\xx).
 \end{align*}
We split the sum according to the annuli $A_m(\yy)$, $m\in\N_0$, which form a disjoint decomposition of $\R^n$. For $\xx\in\frakX_{\ell}\cap A_m(\yy)$ we have $\|\xx-\yy\|_2<(m+1)\delta_{\ell}$ and hence, with \eqref{eq:G},
 \begin{align*}
|\ce_{\ell-1}(\xx)|\leq|\ce_{\ell-1}(\yy)|+|g(\xx)-g(\yy)|\leq\varepsilon_{\ell-1}+\omega_{(m+1)\delta_{\ell}}(\yy).
 \end{align*}
Note that the right-hand side only depends on $\yy$ and $m$, but not on $\xx$. Thus, \eqref{eq:annulus} yields
\begin{align*}
|\ce_{\ell}(\yy)|&\leq|\ce_{\ell-1}(\yy)|+\sum_{m=0}^\infty\sum_{\xx\in\frakX_{\ell}\cap A_m(\yy)}|L_{\xx}^{\frakX_{\ell}}(\yy)|\,|\ce_{\ell-1}(\xx)|\\
&\leq\varepsilon_{\ell-1}+\sum_{m=0}^\infty \left( \sum_{\xx\in\frakX_{\ell} \cap A_m(\yy)} |L_{\xx}^{\frakX_{\ell}}(\yy)| \right) \left( \varepsilon_{\ell-1}+\omega_{(m+1)\delta_{\ell}}(\yy)\right)\\
&\leq\varepsilon_{\ell-1}+\sum_{m=0}^\infty b_m \left( \varepsilon_{\ell-1} + \omega_{(m+1)\delta_{\ell}}(\yy) \right)=(1+\Lambda)\varepsilon_{\ell-1}+\sum_{m=0}^\infty b_m\,\omega_{(m+1)\delta_{\ell}}(\yy).
\end{align*}
Clearly, this bound also holds if $\frakX_{\ell}=\emptyset$. The Cauchy-Schwarz inequality gives
 \begin{align*}
\left( \sum_{m=0}^\infty b_m \omega_{(m+1)\delta_{\ell}}(\yy) \right)^2 \leq \Lambda\sum_{m=0}^\infty b_m\omega_{(m+1)\delta_{\ell}}(\yy)^2.
 \end{align*}
Hence, using $(a+b)^2\leq2a^2+2b^2$ and $\frakD_{\ell}\subseteq X_{\ell}$, we have
 \begin{align*}
\rho_{\ell}^2\leq2(1+\Lambda)^2\varepsilon_{\ell-1}^2h_{\ell}^n\#X_{\ell}+2\Lambda\sum_{m=0}^\infty b_m\,h_{\ell}^n\sum_{\yy\in X_{\ell}}\omega_{(m+1)\delta_{\ell}}(\yy)^2.
 \end{align*}
We apply \cref{lem:oscillation} with $r=(m+1)\delta_{\ell}\geq h_{\ell}$ and use $(m+1)^2\leq(m+1)^{2s}$, which holds since $s\geq1$. This gives
 \begin{align*}
\rho_{\ell}^2\leq2C_\#(1+\Lambda)^2\varepsilon_{\ell-1}^2+2\Lambda C_OB_s\left(\delta_{\ell}^2|g|_{H^1(\R^n)}^2+\delta_{\ell}^{2s}|g|_{H^s(\R^n)}^2\right),
 \end{align*}
which proves the claim.
\end{proof}

The second term in \cref{prop:rho} measures the variation of $E\ce_{\ell-1}$ on balls with radius proportional to $\delta_{\ell}$. By \eqref{eq:normequiv}, it is bounded by $C\|E\ce_{\ell-1}\|_{\Phi_{\ell}}$, but it is not small compared with this norm, so \cref{prop:rho} alone yields no contraction. The smallness comes from the previous level: by \eqref{eq:stab-a} on level $\ell-1$, the $H^s$-seminorm of $\ce_{\ell-1}$ is controlled on the scale $\delta_{\ell-1}$ by $\|E\ce_{\ell-2}\|_{\Phi_{\ell-1}}$. Measured on the finer scale $\delta_{\ell}=\mu\delta_{\ell-1}$, this gives a factor $\mu^s$, and a factor $\mu$ for the $H^1$-term. Hence, $\rho_{\ell}$ is small compared with the residual of the previous level, but not compared with the current one. This is made precise in the next lemma, and it is the reason why the recursion in \cref{thm:main} couples two consecutive levels.

\begin{lemma}\label{lem:coarse}
Let \cref{ass:standing} hold. Then the estimate
 \begin{align*}
\delta_{\ell}\|\ce_{\ell-1}\|_{H^1(\Omega)}+\delta_{\ell}^s\|\ce_{\ell-1}\|_{H^s(\Omega)}\leq K_3\,\mu\left(\|E\ce_{\ell-2}\|_{\Phi_{\ell-1}}+\rho_{\ell-1}\right)
 \end{align*}
holds for all $\ell\geq2$ with $K_3:=C_I\left(c_2^{1/2}(1+C_S)+C_S\right)+c_2^{1/2}$.
\end{lemma}

\begin{proof}
\Cref{lem:stability} on level $\ell-1$ gives
 \begin{align*}
\delta_{\ell-1}^s\|\ce_{\ell-1}\|_{H^s(\Omega)}\leq c_2^{1/2}\|E\ce_{\ell-2}\|_{\Phi_{\ell-1}}\quad\text{and}\quad\|\ce_{\ell-1}\|_{L_2(\Omega)}\leq C_S\left(c_2^{1/2}\|E\ce_{\ell-2}\|_{\Phi_{\ell-1}}+\rho_{\ell-1}\right).
 \end{align*}
Since $\delta_{\ell} \leq \mu\delta_{\ell-1}$, $\mu<1$ and $s\geq1$, the first estimate yields
 \begin{align*}
\delta_{\ell}^s\|\ce_{\ell-1}\|_{H^s(\Omega)}\leq c_2^{1/2}\mu^s\|E\ce_{\ell-2}\|_{\Phi_{\ell-1}}\leq c_2^{1/2}\mu\,\|E\ce_{\ell-2}\|_{\Phi_{\ell-1}}.
 \end{align*}
For the second term, we use \eqref{eq:interpolation} and Young's inequality in the form $A^{1/s}B^{1-1/s}\leq A+B$, which gives
\begin{align*}
\delta_{\ell}\|\ce_{\ell-1}\|_{H^1(\Omega)}&\leq C_I\,\mu\left(\delta_{\ell-1}^s\|\ce_{\ell-1}\|_{H^s(\Omega)}\right)^{1/s}\|\ce_{\ell-1}\|_{L_2(\Omega)}^{1-1/s}\\
&\leq C_I\,\mu\left(c_2^{1/2}\|E\ce_{\ell-2}\|_{\Phi_{\ell-1}}+C_Sc_2^{1/2}\|E\ce_{\ell-2}\|_{\Phi_{\ell-1}}+C_S\rho_{\ell-1}\right).
\end{align*}
Adding both estimates proves the claim.
\end{proof}

Combining \cref{prop:rho}, the boundedness of the extension operator and \cref{lem:coarse}, we obtain with $K:=K_2C_EK_3$
\begin{equation}\label{eq:rho-recursion}
\rho_{\ell}\leq K_1\varepsilon_{\ell-1}+K\mu\left(\|E\ce_{\ell-2}\|_{\Phi_{\ell-1}}+\rho_{\ell-1}\right),\quad \ell\geq2.
\end{equation}
Together with \eqref{eq:recursion}, this is a coupled system of recursions, which we use to prove our main result.

\begin{remark}\label{rem:nu-K}
The constants $C_P$ in \cref{lem:recursion} and $K_3$ in \cref{lem:coarse} are independent of $\mu$ and $\nu$, and so is $C_O$ in \cref{lem:oscillation}. The constants $K_1=(2C_\#)^{1/2}(1+\Lambda)$ and $K_2=(2\Lambda C_OB_s)^{1/2}$ in \cref{prop:rho} are independent of $\mu$, but depend on $\nu$ through $\Lambda$ and $B_s$. By \cref{rem:nu-localization}, we have
\begin{align*}
K_1 \leq C\nu^{2s+n+s(n+1)},\quad \text{and} \quad K=K(\nu)=K_2C_EK_3 \leq C\nu^{p},
\end{align*}
with $ p:=n+s(n+s+3)$, and $C$ depending only on $\Phi$, $n$, $s$ and $c_{qu}$. Note that $K_1$ only multiplies the thresholds, while $K$ enters the contraction factor.
\end{remark}

The mechanism behind the next theorem is the following. At a dropped site the residual of the preceding level is bounded by the threshold. The new local interpolant at this site is a localized average of values of $\ce_{\ell-1}$ at nearby selected sites, and these values differ from $\ce_{\ell-1}(\yy)$ only by the oscillation of $\ce_{\ell-1}$ on the scale $\delta_{\ell}$, which is of order $\mu$. This also explains why $\|\ce_{\ell}\|_{\ell_\infty(\frakD_{\ell})}$ can be much larger than $\varepsilon_{\ell-1}$ while $\rho_{\ell}$ remains small: isolated large values are allowed, a large discrete $\ell_2$-mass is not.

\begin{theorem}\label{thm:main}
Let \cref{ass:standing} hold. Let $f\in H^s(\Omega)$, let $\varepsilon_{\ell}\geq0$ and let $\cf_L$ be computed by \cref{alg:adaptive}. Set
\begin{align*}
r(\mu,\nu):=\widetilde{C}\left(\mu^s+\nu^{-s}\right)+(C_P+2K)\sqrt{\mu}.
\end{align*}
Then there exists a constant $C>0$, independent of $f$, $L$ and the thresholds, but dependent on $ \nu $, such that
\begin{align}\label{eq:main}
\| f - \cf_L\|_{L_2(\Omega)} \leq C \left( r(\mu,\nu)^{L} \|f\|_{H^s(\Omega)} + \mu^{-1/2} \sum_{k=1}^{L-1} r(\mu,\nu)^{L-k} \varepsilon_k \right)
\end{align}
holds for all $L\in\N$. The constants $\widetilde{C}$ and $C_P$ are independent of $\mu$ and $\nu$, and $K=K(\nu)$ is independent of $\mu$. In particular, if $\nu \geq \nu_0:=(3\widetilde{C})^{1/s}$, then $r(\mu,\nu)\leq2/3$ for all
\begin{align}\label{eq:mu0}
\mu \leq \mu_0(\nu) := \min\left\{ (6\widetilde{C})^{-1/s},\,\left(6(C_P+2K(\nu)) \right)^{-2}\right\}.
\end{align}
\end{theorem}

\begin{proof}
For $\ell\in\N$ we define the vector
\begin{align*}
\ww_{\ell}:=\begin{pmatrix}\|E\ce_{\ell-1}\|_{\Phi_{\ell}}\\\rho_{\ell}\end{pmatrix}\in\R^2, \quad \ww_1=\begin{pmatrix}\|Ef\|_{\Phi_1}\\0\end{pmatrix}.
\end{align*}
We can express the recursions \eqref{eq:recursion} and \eqref{eq:rho-recursion} componentwise as
\begin{align*}
\ww_{\ell+1} \leq M\ww_{\ell} + \begin{pmatrix} 0 \\ K_1\varepsilon_{\ell} \end{pmatrix},\quad \text{with} \quad M:=\begin{pmatrix} \alpha & C_P \\ K\mu & K\mu \end{pmatrix}.
\end{align*}
We use the monotone norm $\|\ww\|_*:=\|D^{-1}\ww\|_\infty$ with $D:=\diag(1,\sqrt{\mu})$, for which the associated matrix norm satisfies
\begin{align*}
\|M\|_*=\|D^{-1}MD\|_\infty=\max\left(\alpha+C_P\sqrt\mu,\,K\sqrt\mu+K\mu\right)\leq r(\mu,\nu).
\end{align*}
Iterating the recursion gives
\begin{align*}
\|\ww_L\|_*\leq r(\mu,\nu)^{L-1}\|Ef\|_{\Phi_1}+\mu^{-1/2}K_1\sum_{k=1}^{L-1}r(\mu,\nu)^{L-1-k}\varepsilon_k.
\end{align*}
Next, we apply \eqref{eq:stab-b} with $\ell=L$. Its right-hand side contains exactly the two components of $\ww_L$, and we have $\|E\ce_{L-1}\|_{\Phi_L}\leq\|\ww_L\|_*$ and $\rho_L\leq\sqrt\mu\,\|\ww_L\|_*$. Since $\widetilde{C}\nu^{-s}\leq r(\mu,\nu)$ and $C_P\sqrt\mu\leq r(\mu,\nu)$, this yields
\begin{align*}
\| f - \cf_L \|_{L_2(\Omega)} \leq C_S\left(c_2^{1/2}\nu^{-s}+\sqrt\mu\right)\|\ww_L\|_*\leq C_S\left(\frac{c_2^{1/2}}{\widetilde{C}}+\frac{1}{C_P}\right)r(\mu,\nu)\,\|\ww_L\|_*.
\end{align*}
Together with the bound on $\|\ww_L\|_*$ and $\|Ef\|_{\Phi_1}\leq c_1^{-1/2}C_E\|f\|_{H^s(\Omega)}$ this proves \eqref{eq:main}.

For the second statement note that $\nu\geq\nu_0$ gives $\widetilde{C}\nu^{-s}\leq1/3$, and $\mu\leq\mu_0(\nu)$ gives $\widetilde{C}\mu^s\leq1/6$ and $(C_P+2K)\sqrt\mu\leq1/6$.
\end{proof}

We can now distinguish two different thresholding strategies: Setting a constant threshold for all levels or choosing a strategy where the threshold shrinks geometrically with the number of levels. 

\begin{corollary}\label{cor:thresholds}
With the notation and assumptions of \cref{thm:main} let $r:=r(\mu,\nu)<1$.
\begin{enumerate}
\item If $\varepsilon_k=\varepsilon$ for all $k\in\N$, then
\begin{align*}
\|f-\cf_L\|_{L_2(\Omega)}\leq Cr^{L}\|f\|_{H^s(\Omega)}+C\mu^{-1/2}\frac{\varepsilon}{1-r}.
\end{align*}
\item If $\varepsilon_k=\varepsilon_0\vartheta^k$ with $\vartheta\in(0,1)$, then
\begin{align*}
 \|f-\cf_L\|_{L_2(\Omega)} \leq C \left( r^L \|f\|_{H^s(\Omega)} + \mu^{-1/2} \varepsilon_0 L\max(r,\vartheta)^L\right),
\end{align*}
i.e., the adaptive multilevel method converges linearly.
\end{enumerate}
\end{corollary}

\subsection{Discussion of the Main Result}\label{subsec:discussion}

The results of \cref{subsec:APriori,subsec:APosteriori} warrant detailed discussion. 

\textbf{What is proved.}
\Cref{thm:aposteriori} bounds the error for arbitrary selections of sites in terms of the computable quantities $\rho_k$, and \cref{thm:main} shows that \cref{alg:adaptive} converges if the overlap parameter $\nu$ is fixed first and the refinement parameter is chosen as $\mu\leq\mu_0(\nu)$. 

\textbf{The parameter restriction.}
By \cref{rem:nu-K}, $\mu_0(\nu)\geq c\nu^{-2p}$ for large $\nu$, with constants which are not known explicitly. Hence, the condition $\mu\leq\mu_0(\nu)$ is sufficient but far from necessary: in all our experiments, \cref{alg:adaptive} converges already for $\mu=1/2$. In the coupled regime $\nu=\gamma/\mu$ the term $K(\nu)\sqrt{\mu}\geq c\mu^{1/2-p}$ does not tend to zero, and \cref{thm:main} does not apply. The experiments in \cref{subsec:CoupledRegime}, however, show no deterioration as $\mu\to0$ with fixed $\gamma$.

\textbf{The dependence on $\nu$.}
The growth of $K(\nu)$ stems from the constants of \cref{lem:localization}, which hold uniformly for all subsets $\frakX\subseteq X_{\ell}$. We conjecture that this uniformity cannot be achieved with constants independent of $\nu$. If the selected set contains holes which are larger than $\delta_{\ell}$, the Lagrange functions extrapolate across a distance of up to $\delta_{\ell}=\nu h_{\ell}$, i.e., across about $\nu$ grid spacings, and we expect their Lebesgue constants to grow polynomially in $\nu$. Hence, a polynomial dependence on $\nu$ is inherent in any analysis which is uniform over arbitrary subsets, and only the exponents in \cref{rem:nu-K} are pessimistic.

\textbf{The rate.}
The factor $\sqrt{\mu}$ in $r(\mu,\nu)$ is caused by the structure of the argument rather than by a single estimate. The dropped residual enters the error with a factor of order one, since $\rho_{\ell}$ is part of the $L_2$-error, which is confirmed by the experiments, where $\rho_{\ell}$ closely follows the error. On the other hand, \cref{prop:rho,lem:coarse} bound $\rho_{\ell}$ by the smoothness of $\ce_{\ell-1}$, which is controlled by the previous level. If $\rho_{\ell}$ is bounded by $K\mu^a$ times quantities of level $\ell-1$, the spectral radius of the resulting matrix $M$ behaves like $\sqrt{C_P K \mu^a}$, i.e., the rate is at best $\mu^{a/2}$. In \cref{prop:rho}, the comparison with constants gives $a=1$. If the Lagrange functions on subsets approximately reproduced polynomials, one would obtain $a=s$, but even then only the rate $\mu^{s/2}$, which is still worse than the rate $\alpha\sim\mu^s+\nu^{-s}$ of the classical method. To obtain the rate $\alpha$, one needs a bound of the form
\begin{align*}
\rho_{\ell}\leq q\,\|E\ce_{\ell-1}\|_{\Phi_{\ell}}+C\varepsilon_{\ell-1}\quad\text{with } q\leq C\alpha,
\end{align*}
in which the residual at the dropped sites is controlled by the current level. This is a local a posteriori estimate for the residual near the boundary of the selected region. The check-loop of \cref{thm:loop} enforces such a bound with $ q = 0 $, and the experiments in \cref{sec:numerics} suggest that \cref{alg:adaptive} behaves in the same way: with geometric thresholds it reproduces the reduction of the classical method, also in the coupled regime. We therefore state the following conjecture.

\begin{conjecture}\label{conj:rate}
Let \cref{ass:standing} hold with $\alpha<1$, let $f\in H^s(\Omega)$, and let $\varepsilon_{\ell}=\varepsilon_0\vartheta^{\ell}$ with $\vartheta\in(0,1)$. Then the approximations computed by \cref{alg:adaptive} satisfy
\begin{align*}
\|f-\cf_L\|_{L_2(\Omega)}\leq C\left(\alpha^L\|f\|_{H^s(\Omega)}+\varepsilon_0L\max(\alpha,\vartheta)^L\right),\quad L\in\N,
\end{align*}
with a constant $C>0$ independent of $f$ and $L$.
\end{conjecture}

A proof requires local error estimates for kernel-based interpolation on subsets with large holes, which, to our knowledge, are not available. Another possible route is the convergence theory of adaptive finite element methods, which establishes the contraction of a weighted sum of the error and an error estimator, see \cite{CasconKreuzerNochettoSiebert2008,CarstensenFeischlPagePraetorius2014}. The norm $\|\cdot\|_*$ in the proof of \cref{thm:main} is of the same type. The missing ingredient is a quasi-orthogonality across levels, which is not available since the kernels $\Phi_{\ell}$, and consequently the inner products, change from level to level. This is an open problem even for the classical multilevel method.

\subsection{Convergence of the Check-Loop Variation}\label{subsec:CheckLoopConvergence}

As it turns out, the check-loop variation of \cref{alg:checkloop} allows us to prove convergence without the additional restriction $\mu\leq\mu_0(\nu)$ of \cref{thm:main}, and with the same convergence rate as the non-adaptive method of \cref{sec:multilevel}. Recall that for the classical multilevel method, \eqref{eq:classical} together with \eqref{eq:SamplingInequality} gives
\begin{align*}
\|f-f_L\|_{L_2(\Omega)}\leq C\alpha^L\|f\|_{H^s(\Omega)},\quad L\in\N,
\end{align*}
with $\alpha=\widetilde{C} (\mu^s+\nu^{-s})$. The proof of the next theorem follows immediately from \cref{thm:aposteriori}.

\begin{theorem}\label{thm:loop}
Let \cref{ass:standing} hold. Let $ \tau_{\ell} > 0 $. Let $ \cf_{L} $ be computed by \cref{alg:checkloop}, then there is a constant $C>0$, independent of $f$, $L$, $\mu$, $\nu$ and the
thresholds, such that for $0 \leq t \leq s$
\begin{align}\label{eq:loop-Ht}
\|f-\cf_L\|_{H^t(\Omega)} &\leq C h_L^{-t} \left(\alpha^L\|f\|_{H^s(\Omega)}
  + \sum_{k=1}^{L}\alpha^{L-k}\tau_k \right), 
\end{align}
and 
\begin{align}\label{eq:loop-Linf}
\|f-\cf_L\|_{L_\infty(\Omega)} &\leq C \left(h_L^{-n/2} \left(\alpha^L\|f\|_{H^s(\Omega)}
  + \sum_{k=1}^{L-1}\alpha^{L-k}\tau_k \right) + \tau_L \right) 
\end{align}
hold for all $L\in\N$.
\end{theorem}

\begin{proof}
After termination we have $|\ce_\ell(\xx)|\le\tau_\ell$ for all
$\xx\in X_\ell$, and hence
 \begin{align*}
\eta_\ell\le\tau_\ell \quad\text{and}\quad
\rho_\ell \leq \left(h_\ell^n \#\frakD_\ell \right)^{1/2} \tau_\ell \leq C_\#^{1/2}\tau_\ell,
 \quad  \ell\in\N,
 \end{align*}
 with $ C_{\#} $ from \eqref{eq:Csharp}.
Inserting this into \cref{thm:aposteriori} proves \eqref{eq:loop-Ht} and
\eqref{eq:loop-Linf}.
\end{proof}

\begin{corollary}\label{cor:checkloop-rate}
Let the assumptions of \cref{thm:loop} hold with $\alpha<1$ and $\tau_{\ell}=\tau_0\vartheta^{\ell}$, $\vartheta\in(0,1)$. Then
\begin{align*}
\|f-\cf_L\|_{L_2(\Omega)}\leq
\begin{cases}
C\left(\|f\|_{H^s(\Omega)}+\dfrac{\tau_0\vartheta}{\alpha-\vartheta}\right)\alpha^L, & \vartheta<\alpha,\\[2mm]
C\left(\|f\|_{H^s(\Omega)}+\tau_0L\right)\alpha^L, & \vartheta=\alpha,\\[1mm]
C\left(\alpha^L\|f\|_{H^s(\Omega)}+\dfrac{\tau_0\vartheta}{\vartheta-\alpha}\,\vartheta^L\right), & \vartheta>\alpha,
\end{cases}
\end{align*}
with the constant $C$ of \cref{thm:loop}, which is independent of $f$, $L$, $\mu$ and $\nu$. In particular, for $\vartheta<\alpha$, \cref{alg:checkloop} converges with exactly the rate $\alpha^L$ of the classical multilevel method, and this also holds in the coupled regime $\nu=\gamma/\mu$, in which $\alpha=\widetilde{C}(1+\gamma^{-s})\mu^s$.
\end{corollary}

\begin{proof}
By \eqref{eq:loop-Ht} with $t=0$, it remains to bound $\sum_{k=1}^{L}\alpha^{L-k}\tau_0\vartheta^k=\tau_0\alpha^L\sum_{k=1}^{L}(\vartheta/\alpha)^k$. For $\vartheta<\alpha$, the geometric series is bounded by $\vartheta/(\alpha-\vartheta)$. For $\vartheta=\alpha$, the sum equals $L$. For $\vartheta>\alpha$, we write it as $\tau_0\vartheta^L\sum_{j=0}^{L-1}(\alpha/\vartheta)^j\leq\tau_0\vartheta^L\,\vartheta/(\vartheta-\alpha)$.
\end{proof}

\section{Extensions}\label{sec:extensions}

The analysis of \cref{sec:adaptive} only uses three properties of the local approximation operator on level $\ell$: it is stable in the native space norm, it controls the residual on the selected sites $\frakX_{\ell}$, and its values at the dropped sites are localized averages of the data on $\frakX_{\ell}$. In this section we verify these properties for two variants of the multilevel method which have been analyzed for the classical method in \cite{Wendland2010}: penalized least-squares approximation and interpolation of target functions which are rougher than the native space of the kernel.

\subsection{Penalized Least-Squares Approximation}

For $\lambda > 0$, a nonempty set $ \frakX \subseteq X_{\ell}$ and $g \in C(\Omega)$ we define the penalized least-squares approximation as
\begin{equation}\label{eq:pls}
Q_{\frakX,\Phi_{\ell}}^\lambda g := \argmin_{s\in \cN_{\Phi_{\ell}}(\R^n)}\left\{\sum_{\xx \in \frakX}| g(\xx) -s(\xx)|^2 + \lambda \|s\|_{\Phi_{\ell}}^2\right\},
\end{equation}
and set $Q_{\emptyset,\Phi_{\ell}}^\lambda g := 0$. It is well-known that $Q_{\frakX,\Phi_{\ell}}^\lambda g = \sum_{\xx \in \frakX} \beta_{\xx}\Phi_{\ell}(\cdot-\xx)\in W_{\ell}(\frakX)$, where the coefficient vector $ \bbeta $ is the solution of $(M_{\frakX}+\lambda I) \bbeta = g(\frakX)$, see, e.g., \cite{Wendland2004}. The adaptive penalized least-squares multilevel method is \cref{alg:adaptive} with the local approximation $\cs_{\ell} = Q_{\frakX_{\ell},\Phi_{\ell}}^{\lambda_{\ell}} (\ce_{\ell-1}) $. We use the same notation $\cs_{\ell}$, $\ce_{\ell}$, $\cf_{\ell}$, $\frakD_{\ell}$ and $\rho_{\ell}$ as before. Note that $\ce_{\ell}$ no longer vanishes on $\frakX_{\ell}$.

The next lemma collects properties needed to prove the convergence result, for a proof see, e.g., \cite{WendlandRieger2005}.

\begin{lemma}\label{lem:pls}
Let $\ell \in \N$, $ \lambda > 0$, $\frakX \subseteq X_{\ell}$  and $g \in \cN_{\Phi_{\ell}}(\R^n) $. Then $s := Q_{\frakX,\Phi_{\ell}}^\lambda g$ satisfies
\begin{align}\label{eq:pls-a}
\|g - s \|_{\Phi_{\ell}} &\leq 2 \| g \|_{\Phi_{\ell}},
\end{align}
and
\begin{align}\label{eq:pls-b}
\sum_{\xx\in\frakX} | g(\xx) - s(\xx) |^2 &\leq \lambda \| g \|_{\Phi_{\ell}}^2.
\end{align}
\end{lemma}

The values of $\cs_{\ell}$ at the dropped sites can again be written as localized averages. For a nonempty set $ \frakX \subseteq X_{\ell}$, $ \lambda > 0 $ and $\xx \in \frakX $ we define
\begin{align*}
L_{\xx}^{\frakX,\lambda}(\yy) := \sum_{\zz\in\frakX} \left((M_{\frakX} + \lambda I)^{-1} \right)_{\xx,\zz} \Phi_{\ell} (\yy-\zz), \quad \yy\in\R^n,
\end{align*}
such that $ Q_{\frakX,\Phi_{\ell}}^\lambda g(\yy) = \sum_{\xx\in\frakX} g(\xx) L_{\xx}^{\frakX,\lambda}(\yy) $. These functions do not satisfy the Lagrange condition, however, this has not been used in the proof of \cref{prop:rho}.

\begin{lemma}\label{lem:localization-pls}
Let \cref{ass:standing} hold. Then the estimate of \cref{lem:localization} holds for $L_{\xx}^{\frakX,\lambda}$ with the same constants $C_L$ and $\theta$, for every $\lambda > 0$.
\end{lemma}

\begin{proof}
The matrix $M_{\frakX}+\lambda I$ has the same sparsity pattern as $M_{\frakX}$ and is therefore $R$-banded with the same $R$. Its eigenvalues are those of $M_{\frakX}$ shifted by $\lambda$. Hence, 
\begin{align*}
\|(M_{\frakX}+\lambda I)^{-1}\|_2\leq\|M_{\frakX}^{-1}\|_2
\end{align*}
and
\begin{align*}
\cond_2(M_{\frakX}+\lambda I) = \frac{\lambda_{\max}(M_{\frakX})+\lambda}{\lambda_{\min}(M_{\frakX})+\lambda}\leq\cond_2(M_{\frakX})\leq\kappa_\nu.
\end{align*}
Since the proof of \cref{lem:localization} uses the matrix only through upper bounds on these two quantities, it applies verbatim.
\end{proof}

We choose the smoothing parameters such that
\begin{align}\label{eq:lambda}
\lambda_{\ell}\leq c_{\lambda}h_{\ell}^{-n}\nu^{-2s}, \quad \ell\in\N,
\end{align}
with a fixed $c_{\lambda}>0$. The factor $h_{\ell}^{-n}$ reflects that the data term in \eqref{eq:pls} is a sum over $\#X_{\ell}\sim h_{\ell}^{-n}$ sites. For the normalized data term $h_{\ell}^n\sum_{\xx\in\frakX} |g(\xx)-s(\xx)|^2$, \eqref{eq:lambda} becomes $\lambda_{\ell}\leq c_{\lambda}\nu^{-2s}$, which is the condition of \cite[Theorem~3]{Wendland2010}. Since we use the $\ell_2$-version \eqref{eq:sampling-l2} of the sampling inequality instead of the $\ell_\infty$-version, \eqref{eq:lambda} is weaker than the condition in \cite{Wendland2010}. 

In terms of the three properties listed at the beginning of this section, \cref{lem:pls} gives the stability in the native space norm, \eqref{eq:pls-b} together with \eqref{eq:lambda} controls the residual on the selected sites, and \cref{lem:localization-pls} gives the localization.

\begin{theorem}\label{thm:pls}
Let \cref{ass:standing} hold and let the smoothing parameters satisfy \eqref{eq:lambda}. Let $\varepsilon_{\ell} \geq 0$, and let $\cf_L$ be computed by the adaptive penalized least-squares multilevel method. Then all results of \cref{sec:adaptive} hold with $\alpha$ replaced by
\begin{align*}
\alpha_\lambda := \widetilde{C} \left( 2 \mu^s + \left(2+(c_\lambda/c_2)^{1/2}\right)\nu^{-s}\right),
\end{align*}
and with constants $ K_1 $, $ K_2 $ as before and modified constants $K_3$, $K$, which are independent of $\mu$ and depend on $\nu$ as in \cref{rem:nu-K}. In particular, \eqref{eq:main} holds with
\begin{align*}
r_\lambda(\mu,\nu):=\widetilde{C}\left(2 \mu^s+\left( 2 +(c_\lambda/c_2)^{1/2}\right)\nu^{-s}\right)+(C_P + 2K)\sqrt\mu
\end{align*}
in place of $r(\mu,\nu)$, and \cref{thm:loop} holds if the check-loop is performed on $\frakD_{\ell}$.
\end{theorem}

\begin{proof}
We go through the results of \cref{sec:adaptive} and point out the differences. Estimate \eqref{eq:stab-a} follows as before, where we use \eqref{eq:pls-a} instead of the projection property of the interpolant. For \eqref{eq:stab-b}, we apply \eqref{eq:sampling-l2} on $X_{\ell}$ and use $\|\ce_{\ell}\|_{\ell_2(X_{\ell})} \leq \|\ce_{\ell}\|_{\ell_2(\frakX_{\ell})} + \|\ce_{\ell}\|_{\ell_2(\frakD_{\ell})}$. With \eqref{eq:pls-b} applied to $E\ce_{\ell-1}$, which coincides with $\ce_{\ell-1}$ on $\frakX_{\ell}$, and \eqref{eq:lambda} we have
\begin{align*}
h_{\ell}^{n/2}\|\ce_{\ell}\|_{\ell_2(\frakX_{\ell})}\leq\left(h_{\ell}^n\lambda_{\ell}\right)^{1/2}\|E\ce_{\ell-1}\|_{\Phi_{\ell}}\leq c_{\lambda}^{1/2}\nu^{-s}\|E\ce_{\ell-1}\|_{\Phi_{\ell}}.
\end{align*}
Hence, \eqref{eq:stab-b} becomes
\begin{align*}
\|\ce_{\ell}\|_{L_2(\Omega)} \leq  C_S\left(\left(2c_2^{1/2}+c_{\lambda}^{1/2}\right)\nu^{-s}\|E\ce_{\ell-1}\|_{\Phi_{\ell}}+\rho_{\ell}\right),
\end{align*}
and the proof of \cref{lem:recursion} yields \eqref{eq:recursion} with $\alpha_{\lambda}$ instead of $\alpha$ and the same $C_P$. In the proof of \cref{thm:aposteriori}, the data terms on $X_L$ now also contain the residual on $\frakX_L$. By the estimate above and $\|\ce_L\|_{\ell_\infty(\frakX_L)}\leq\|\ce_L\|_{\ell_2(\frakX_L)}$, these parts are bounded by $c_{\lambda}^{1/2}h_L^{-t}\nu^{-s}\|E\ce_{L-1}\|_{\Phi_L}$ and $c_{\lambda}^{1/2}h_L^{-n/2}\nu^{-s}\|E\ce_{L-1}\|_{\Phi_L}$, respectively, and are absorbed into the first term by \eqref{eq:aposteriori-aux}. In the proof of \cref{prop:rho} we use the representation $\cs_{\ell}(\yy)=\sum_{\xx\in\frakX_{\ell}}L_{\xx}^{\frakX_{\ell},\lambda_{\ell}}(\yy)\ce_{\ell-1}(\xx)$ and \cref{lem:localization-pls} instead of \cref{lem:localization}, which gives the same $\Lambda$ and $B_s$, and hence the same $K_1$ and $K_2$. \Cref{lem:coarse} follows from the modified version of \cref{lem:stability} on level $\ell-1$, with $c_2^{1/2}$ replaced by $2c_2^{1/2} $ in the bound for the $ H^s $-seminorm and $ 2 c_2^{1/2} + c_{\lambda}^{1/2} $ in the bound for the $L_2$-norm, which only changes $K_3$. Finally, the proof of \cref{thm:main} only uses \eqref{eq:recursion}, \eqref{eq:rho-recursion} and \eqref{eq:stab-b}. For \cref{thm:loop} note that $\rho_{\ell}$ and $\eta_{\ell}$ only involve $\frakD_{\ell}$, such that it suffices to check the residual on $\frakD_{\ell}$.
\end{proof}

\subsection{Escaping the Native Space}\label{subsec:escaping}

Next, we consider target functions $f\in H^\sigma(\Omega)$ with $n/2<\sigma\leq s$, i.e., functions that are rougher than the native space of $\Phi$. The method itself is unchanged, i.e., we use \cref{alg:adaptive} with the kernels $\Phi_{\ell}$. For the analysis, however, we follow \cite[Section~5]{Wendland2010} and measure the residuals in the native space norms of a second kernel $\Psi:\R^n\to\R$, which is a reproducing kernel of $H^\sigma(\R^n)$, i.e., its Fourier transform satisfies \eqref{eq:fourier} with $\sigma$ instead of $s$. Without loss of generality we use the same equivalence constants $c_1,c_2$ as for $\Phi$. We set $\Psi_{\ell}:=\delta_{\ell}^{-n}\Psi(\cdot/\delta_{\ell})$. Then \eqref{eq:normequiv} holds for $\Psi_{\ell}$ with $\sigma$ instead of $s$. The key ingredient is the following stability estimate, which again holds uniformly for all subsets of $X_{\ell}$.

\begin{lemma}\label{lem:escape}
Let \cref{ass:standing} hold and let $n/2 < \sigma \leq s$. Then there exists a constant $C_\sigma > 0$, depending only on $\Phi$, $\Psi$, $n$, $s$, $\sigma$, $\nu$ and $c_{qu}$, such that
\begin{align*}
\| g - I_{\frakX,\Phi_{\ell}} g \|_{\Psi_{\ell}} \leq C_\sigma \| g \|_{\Psi_{\ell}}
\end{align*}
holds for all $ \ell \in \N$, all $ \frakX \subseteq X_{\ell}$ and all $ g \in H^\sigma(\R^n)$.
\end{lemma}

\begin{proof}
For $\frakX=\emptyset$ there is nothing to show, so let $\frakX\neq\emptyset$. We set $\delta:=\delta_{\ell}$ and $\widetilde{g} := g (\delta\,\cdot)$. A direct computation gives
\begin{align*}
\| \widetilde{g} \|_{H^\sigma(\R^n)}^2 = \delta^{-n} \int_{\R^n}|\widehat{g}(\xxi)|^2 \left( 1 + \delta^2 \|\xxi\|_2^2 \right)^\sigma \ d\xxi,
\end{align*}
and with \eqref{eq:fourier} for $\Psi$ we obtain $c_2^{-1} \delta^n \|\widetilde{g} \|_{H^\sigma(\R^n)}^2 \leq \| g \|_{\Psi_{\ell}}^2 \leq c_1^{-1} \delta^n \| \widetilde{g} \|_{H^\sigma(\R^n)}^2 $. The rescaled set $ \widetilde{\frakX} := \delta^{-1} \frakX $ has separation radius $q_{\widetilde{\frakX}} \geq q_{\ell}/\delta \geq (\nu c_{qu})^{-1}$, since the sites of $ \frakX \subseteq X_{\ell}$ are at least $ 2q_{\ell} $ apart. Hence, $ \widetilde{\frakX} $ is separated with the parameter $ q_* := \min(q_{\widetilde{\frakX}},1) \in [(\nu c_{qu})^{-1},1]$. By \cite{NarcowichWardWendland2006} there exists a band-limited function $ \widetilde{g}_B$ with $ \supp \widehat{\widetilde{g}_B} \subseteq B_{c_B/q_*}(\00)\subseteq B_{c_B\nu c_{qu}}(\00)$, which coincides with $\widetilde{g} $ on $ \widetilde{\frakX} $ and satisfies $ \| \widetilde{g}_B \|_{H^\sigma(\R^n)} \leq C\| \widetilde{g} \|_{H^\sigma(\R^n)} $, where $c_B$ and $C$ only depend on $n$ and $\sigma$. Note that this construction only uses the separation of $ \widetilde{\frakX} $ and not its fill distance, which is the reason why the estimate holds uniformly for all subsets $ \frakX \subseteq X_{\ell} $. We set $ g_B := \widetilde{g}_B (\cdot/\delta) $. Then $ g_B = g $ on $ \frakX $ and therefore $ I_{\frakX,\Phi_{\ell}} g = I_{\frakX,\Phi_{\ell}} g_B$. Since $\Psi$ satisfies \eqref{eq:fourier} with $\sigma \leq s$, we have $\|u\|_{\Psi_{\ell}}\leq(c_2/c_1)^{1/2}\|u\|_{\Phi_{\ell}}$ for $u\in H^s(\R^n)$. Together with the norm minimality of the interpolant this yields
\begin{align*}
\| I_{\frakX,\Phi_{\ell}} g \|_{\Psi_{\ell}} \leq \left(\frac{c_2}{c_1}\right)^{1/2} \| I_{\frakX,\Phi_{\ell}} g_B \|_{\Phi_{\ell}} \leq \left(\frac{c_2}{c_1}\right)^{1/2} \| g_B \|_{\Phi_{\ell}}.
\end{align*}
Since $ \widehat{g_B} $ is supported in $ \|\xxi\|_2 \leq c_B\nu c_{qu}/\delta$, we have $ ( 1 + \delta^2 \|\xxi\|_2^2)^s \leq (1+c_B^2\nu^2c_{qu}^2)^{s-\sigma} (1 + \delta^2 \|\xxi\|_2^2)^\sigma $ on the support. Hence,
\begin{align*}
\| g_B \|_{\Phi_{\ell}} \leq \left(\frac{c_2}{c_1}\right)^{1/2} \left(1 + c_B^2\nu^2c_{qu}^2\right)^{(s-\sigma)/2} \| g_B \|_{\Psi_{\ell}} \leq C\left(1 + c_B^2\nu^2c_{qu}^2\right)^{(s-\sigma)/2} \| g \|_{\Psi_{\ell}},
\end{align*}
where we used the scaling relation and the stability of $\widetilde{g}_B$ in the last step. The claim follows with the triangle inequality.
\end{proof}

\begin{theorem}\label{thm:escape}
Let \cref{ass:standing} hold, let $n/2<\sigma\leq s$ with $\sigma \geq 1$, and let $ f \in H^\sigma(\Omega)$. Let $\varepsilon_{\ell} \geq 0$ and let $\cf_L$ be computed by \cref{alg:adaptive}. Let $\widetilde{C}$ and $C_P$ denote the constants of \cref{lem:recursion} with $\sigma$ instead of $s$, and set
\begin{align*}
\alpha_{\sigma} := \widetilde{C} C_\sigma \left(\mu^\sigma+\nu^{-\sigma} \right),\quad \text{and} \quad r_{\sigma}(\mu,\nu) := \alpha_{\sigma} + \left(C_P + 2 K_\sigma\right)\sqrt{\mu},
\end{align*}
with $C_\sigma$ from \cref{lem:escape} and a constant $K_\sigma=K_\sigma(\nu)$, which is independent of $\mu$. Then there exists a constant $C>0$, independent of $f$, $L$ and the thresholds, such that
\begin{align*}
\| f - \cf_L \|_{L_2(\Omega)} \leq C \left( r_\sigma(\mu,\nu)^L \|f\|_{H^\sigma(\Omega)} + \mu^{-1/2} \sum_{k=1}^{L-1}r_\sigma(\mu,\nu)^{L-k}\varepsilon_k\right)
\end{align*}
holds for all $L\in\N$. Moreover, \cref{thm:aposteriori,thm:loop} hold for $0\leq t\leq\sigma$ with $\alpha_\sigma$ instead of $\alpha$ and $\|f\|_{H^\sigma(\Omega)}$ instead of $\|f\|_{H^s(\Omega)}$, with constants independent of $\mu$ and $\nu$.
\end{theorem}

\begin{proof}
We replace $s$ by $\sigma$ and $\|\cdot\|_{\Phi_{\ell}}$ by $\|\cdot\|_{\Psi_{\ell}}$ in all estimates of \cref{sec:convergence} which concern the residuals. Since $\cs_{\ell}=I_{\frakX_{\ell},\Phi_{\ell}}(E\ce_{\ell-1})$, the function $E\ce_{\ell-1}-\cs_{\ell}$ extends $\ce_{\ell}$, and \eqref{eq:normequiv} for $\Psi_{\ell}$ together with \cref{lem:escape} gives
\begin{align*}
\|\ce_{\ell}\|_{H^\sigma(\Omega)}\leq c_2^{1/2}\delta_{\ell}^{-\sigma}\|E\ce_{\ell-1}-I_{\frakX_{\ell},\Phi_{\ell}}(E\ce_{\ell-1})\|_{\Psi_{\ell}}\leq c_2^{1/2}C_\sigma\delta_{\ell}^{-\sigma}\|E\ce_{\ell-1}\|_{\Psi_{\ell}}.
\end{align*}
The sampling inequality \eqref{eq:sampling-l2} also holds with $\sigma$ instead of $s$, which gives
\begin{align*}
\|\ce_{\ell}\|_{L_2(\Omega)}\leq C_S\left(c_2^{1/2}C_\sigma\nu^{-\sigma}\|E\ce_{\ell-1}\|_{\Psi_{\ell}}+\rho_{\ell}\right).
\end{align*}
These are \eqref{eq:stab-a} and \eqref{eq:stab-b} with the additional factor $C_\sigma$ in front of $\|E\ce_{\ell-1}\|_{\Psi_{\ell}}$. Hence, the proof of \cref{lem:recursion} yields \eqref{eq:recursion} with $\alpha_\sigma$ and the same $C_P$. Since $\alpha_\sigma\geq \widetilde{C}C_\sigma\nu^{-\sigma}$, the proofs of \cref{thm:aposteriori,thm:loop} carry over with constants independent of $\mu$ and $\nu$. \Cref{lem:localization} concerns the kernel matrices of $\Phi_{\ell}$ and remains unchanged, and hence so do $\Lambda$ and $K_1$. \Cref{lem:oscillation} holds with $\sigma$ instead of $s$, since we only used $s>n/2$ and $s\geq1$, and so does \cref{prop:rho}, with $B_\sigma$ instead of $B_s$ in $K_2$. In the proof of \cref{lem:coarse} the two estimates above on level $\ell-1$ and \eqref{eq:interpolation} with $\sigma$ instead of $s$ give the constant
\begin{align*}
K_{3,\sigma}:=C_I\left(c_2^{1/2}C_\sigma(1+C_S)+C_S\right)+c_2^{1/2}C_\sigma,
\end{align*}
and we set $K_\sigma:=K_2C_EK_{3,\sigma}$. The proof of \cref{thm:main} then carries over verbatim, where in the last step we use $C_S(c_2^{1/2}C_\sigma\nu^{-\sigma}+\sqrt\mu)\leq C_S(c_2^{1/2}\widetilde{C}^{-1} + C_P^{-1})r_\sigma(\mu,\nu)$. Finally, $\|Ef\|_{\Psi_1}\leq c_1^{-1/2}C_E\|f\|_{H^\sigma(\Omega)}$.
\end{proof}

\begin{remark}
\begin{enumerate}
   \item The proof of \cref{lem:escape} gives $C_\sigma\leq C\nu^{s-\sigma}$ for $\nu c_{qu}\geq1$. Hence, the term $\widetilde{C}C_\sigma\nu^{-\sigma}$ in $\alpha_\sigma$ is bounded by $C\nu^{s-2\sigma}$, and it becomes small for large $\nu$ only if $\sigma>s/2$. As in \cref{rem:nu-K}, the constant $K_\sigma(\nu)$ grows polynomially in $\nu$.
\item Both extensions can be combined, i.e., penalized least-squares approximation of target functions in $H^\sigma(\Omega)$. Since \cref{lem:pls} has to be replaced by a stability estimate in the $\Psi_{\ell}$-norms, we do not pursue this here.
\end{enumerate}
\end{remark}

\section{Numerical Examples}\label{sec:numerics}

We illustrate our findings with three examples in two dimensions: Franke's function on the unit square, a corner singularity on an L-shaped domain, and the composite test function of \cite{AvesaniIskeMultererWenzel2026}, which exhibits several distinct features, such as singularities along curves.

Throughout, the sets of sites $X_\ell$ are nested uniform grids on $\Omega$, where $X_{\ell+1}$ is obtained from $X_\ell$ by halving the distance between neighboring sites. Hence, $h_{\ell+1}=h_\ell/2$, i.e., $\mu=1/2$. In the implementation, we use $\widetilde h_\ell:=(|\Omega|/N_\ell)^{1/2}$ with $N_\ell:=\#X_\ell$ and set $\delta_\ell=\nu\widetilde h_\ell$. Since $\widetilde h_\ell$ is proportional to $h_\ell$ up to a factor which tends to $\sqrt{2}$ as $\ell$ grows, this means that the overlap parameter in the theory is $ \sqrt{2} $ times the $\nu$ used in the experiments.

Unless stated otherwise, we use $\nu=3$ and $L=8$ levels. We use Wendland's $ C^0 $- and $C^2 $-functions 
\begin{align*}
\phi_{2,0}(r)=(1-r)_+^2, \quad \text{and} \quad \phi_{2,1}(r)=(1-r)_+^4(4r+1),
\end{align*}
which are reproducing kernels of $H^{3/2}(\R^2)$ and $H^{5/2}(\R^2)$, respectively. The error is measured on $ X_{L+1} $, the sites of level $L+1$. We report the root mean square error 
\begin{align*}
\text{rms error} := \left( \frac{1}{N_{L+1}}\sum_{i=1}^{N_{L+1}}|f(\xx_i)-\cf_L(\xx_i)|^2 \right)^{\frac{1}{2}} \approx \frac{1}{|\Omega|^{1/2}}\|f-\cf_L\|_{L_2(\Omega)},   
\end{align*}
and the maximum error. Since a quarter of these points are sites of $X_L$, at which the classical multilevel interpolant is exact, the root mean square error of the classical method is slightly underestimated.

We compare the classical multilevel method, which uses all sites on every level, with \cref{alg:adaptive} for two choices of thresholds: constant thresholds $\varepsilon_{\ell}=\varepsilon$ and geometric thresholds $\varepsilon_{\ell}=\varepsilon_0\vartheta^{\ell}$ where we set $\varepsilon_0=10^{-2}$, see \cref{cor:thresholds}. Moreover, we test the variant with check-loop of \cref{alg:checkloop} with $\tau_{\ell}=\varepsilon_{\ell}/2$. Besides the errors, we report the number of selected sites $\#\frakX_{\ell}$, the total number of centers $\sum_{\ell=1}^L\#\frakX_{\ell}$, and the quantities $\rho_{\ell}$ and $\eta_{\ell}:=\|\ce_{\ell}\|_{\ell_\infty(\frakD_{\ell})}$. To compare the efficiency of the methods, we define the \emph{saving} of an adaptive run as the number of centers which the classical method needs to reach the same final root mean square error, divided by the number of centers used by the adaptive run. Here, the error of the classical method is interpolated linearly in a log-log scale between its levels. A saving larger than one means that the adaptive method needs fewer centers for the same accuracy.

\subsection{Franke's Function}

Our first example is Franke's function \cite{Franke1982} on $\Omega=[0,1]^2$, which is a sum of four exponential functions and hence smooth, in particular $f\in H^{5/2}(\Omega)$. Choosing the kernel $\phi_{2,1}$ means that the target function is inside the native space and \cref{thm:main} applies. The sets of sites have $N_{\ell}=(2^{\ell+1}+1)^2$ points, i.e., $N_1=25$ and $N_8=263\,169$.

\Cref{tab:franke} shows the results of the classical method, of \cref{alg:adaptive} with the constant threshold $\varepsilon=10^{-4}$, and with the geometric thresholds for $\vartheta=0.3$. The classical method reduces the error by a factor of about $0.23$ per level. With the constant threshold, the number of selected sites remains of the order of $10^3$ on the finer levels, and the error stagnates at about $0.3\,\varepsilon$, as predicted by \cref{cor:thresholds}. The quantity $\rho_{\ell}$ closely follows the root mean square error on all levels, while $\eta_{\ell}$ exceeds the threshold by up to two orders of magnitude. This confirms that the $\ell_2$-quantity $\rho_{\ell}$, and not the residual at the dropped sites in the maximum norm, is the right measure for the analysis. With geometric thresholds, the error keeps decreasing at the classical rate, which is better than the rate guaranteed by \cref{cor:thresholds}, see \cref{conj:rate}. Since the threshold decreases at the same rate as the residual, however, the number of selected sites grows with $N_\ell$. For a prescribed target accuracy $ \varepsilon_{target} $ of the approximation, a constant threshold of $ \approx 3 \varepsilon_{target} $  seems to be the most economical.

\begin{table}[ht]
\centering
\small
\begin{tabular}{rr c rccc rc}
\toprule
& & classical & \multicolumn{4}{c}{constant, $\varepsilon=10^{-4}$} & \multicolumn{2}{c}{geometric, $\vartheta=0.3$}\\
\cmidrule(lr){3-3}\cmidrule(lr){4-7}\cmidrule(lr){8-9}
$\ell$ & $N_{\ell}$ & rms error & $\#\frakX_{\ell}$ & rms error & $\rho_{\ell}$ & $\eta_{\ell}$ & $\#\frakX_{\ell}$ & rms error\\
\midrule
1 & 25 & $2.8\cdot10^{-2}$ & 25 & $2.8\cdot10^{-2}$ & 0 & 0 & 25 & $2.8\cdot10^{-2}$\\
2 & 81 & $7.1\cdot10^{-3}$ & 56 & $8.4\cdot10^{-3}$ & $5.6\cdot10^{-3}$ & $2.8\cdot10^{-2}$ & 48 & $8.9\cdot10^{-3}$\\
3 & 289 & $1.2\cdot10^{-3}$ & 228 & $1.3\cdot10^{-3}$ & $8.4\cdot10^{-4}$ & $6.8\cdot10^{-3}$ & 200 & $1.6\cdot10^{-3}$\\
4 & 1\,089 & $2.7\cdot10^{-4}$ & 650 & $2.9\cdot10^{-4}$ & $3.1\cdot10^{-4}$ & $4.7\cdot10^{-3}$ & 548 & $3.4\cdot10^{-4}$\\
5 & 4\,225 & $6.2\cdot10^{-5}$ & 902 & $9.7\cdot10^{-5}$ & $8.8\cdot10^{-5}$ & $1.0\cdot10^{-3}$ & 1\,567 & $9.5\cdot10^{-5}$\\
6 & 16\,641 & $1.4\cdot10^{-5}$ & 1\,225 & $3.6\cdot10^{-5}$ & $4.0\cdot10^{-5}$ & $7.2\cdot10^{-4}$ & 4\,495 & $2.3\cdot10^{-5}$\\
7 & 66\,049 & $3.2\cdot10^{-6}$ & 1\,132 & $3.0\cdot10^{-5}$ & $3.0\cdot10^{-5}$ & $2.5\cdot10^{-4}$ & 17\,769 & $5.7\cdot10^{-6}$\\
8 & 263\,169 & $7.1\cdot10^{-7}$ & 1\,294 & $2.8\cdot10^{-5}$ & $2.8\cdot10^{-5}$ & $1.5\cdot10^{-4}$ & 63\,255 & $1.4\cdot10^{-6}$\\
\bottomrule
\end{tabular}
\caption{Franke’s function with $\phi_{2,1}$, $\mu=1/2$ and $\nu=3$. For each level $\ell$, the number of sites $N_\ell$, the number of selected sites $\#\frakX_\ell$ and the root mean square error of $ f - \cf_\ell$ are given for the classical method, for \cref{alg:adaptive} with the constant threshold $\varepsilon=10^{-4}$, and for the geometric thresholds $\varepsilon_\ell=10^{-2}\cdot0.3^\ell$. For the constant threshold, $\rho_\ell$ and $\eta_\ell $ measure the residual at the dropped sites.}\label{tab:franke}
\end{table}

\Cref{tab:thresholds} compares several thresholds. Here, we also report the mean of the numerical contraction factor per level over the last three levels. For the geometric thresholds, \cref{conj:rate} predicts a reduction of $\max(\alpha,\vartheta)$ per level, whereas \cref{cor:thresholds} only guarantees $\max(r,\vartheta)$. Indeed, for $\vartheta\in\{0.1,0.2\}$, which is below the observed classical reduction of about $0.23$, the adaptive method converges with the classical rate, while for $\vartheta=0.5$ it converges with the rate $\vartheta$. The most efficient choice is $\vartheta$ close to the classical reduction, which yields a saving of about two. For constant thresholds, the saving increases as the threshold decreases and reaches $6.3$ for $\varepsilon=10^{-5}$. The check-loop reduces the maximum error and yields a saving of $3.5$, at the price of two to three additional solves per level.

\begin{table}[ht]
\centering
\small
\begin{tabular}{l rccc c}
\toprule
method & centers & rms error & max error & reduction & saving\\
\midrule
classical & 351\,568 & $7.1\cdot10^{-7}$ & $5.5\cdot10^{-5}$ & 0.23 & --\\
\midrule
constant, $\varepsilon=10^{-3}$ & 1\,030 & $2.6\cdot10^{-4}$ & $9.8\cdot10^{-4}$ & 0.96 & 1.47\\
constant, $\varepsilon=10^{-4}$ & 5\,512 & $2.8\cdot10^{-5}$ & $1.5\cdot10^{-4}$ & 0.66 & 2.17\\
constant, $\varepsilon=10^{-5}$ & 15\,323 & $2.9\cdot10^{-6}$ & $1.1\cdot10^{-4}$ & 0.35 & 6.29\\
\midrule
geometric, $\vartheta=0.1$ & 249\,161 & $8.8\cdot10^{-7}$ & $1.1\cdot10^{-4}$ & 0.23 & 1.16\\
geometric, $\vartheta=0.2$ & 122\,633 & $9.2\cdot10^{-7}$ & $1.1\cdot10^{-4}$ & 0.23 & 2.26\\
geometric, $\vartheta=0.3$ & 87\,907 & $1.4\cdot10^{-6}$ & $1.1\cdot10^{-4}$ & 0.25 & 2.09\\
geometric, $\vartheta=0.5$ & 62\,170 & $2.6\cdot10^{-5}$ & $1.6\cdot10^{-4}$ & 0.50 & 0.21\\
\midrule
check-loop, $\varepsilon=10^{-4}$, $\tau=5\cdot10^{-5}$ & 5\,889 & $1.5\cdot10^{-5}$ & $7.1\cdot10^{-5}$ & 0.62 & 3.52\\
\bottomrule
\end{tabular}
\caption{Franke's function with $\phi_{2,1}$, $\mu=1/2$, $\nu=3$ and $L=8$: the classical method and \cref{alg:adaptive} with constant thresholds $\varepsilon_\ell=\varepsilon$, geometric thresholds $\varepsilon_\ell=10^{-2}\vartheta^\ell$, and the check-loop with $\tau_\ell=\varepsilon_\ell/2$. Given are the total number of centers $\sum_{\ell=1}^L\#\frakX_\ell$, the root mean square and maximum error of $\cf_L$, the mean reduction of the root mean square error per level over the last three levels, and the saving, i.e., the number of centers the classical method needs to reach the same root mean square error divided by the number of centers used.}\label{tab:thresholds}
\end{table}

\begin{figure}[ht]
\centering
\includegraphics[width=\textwidth]{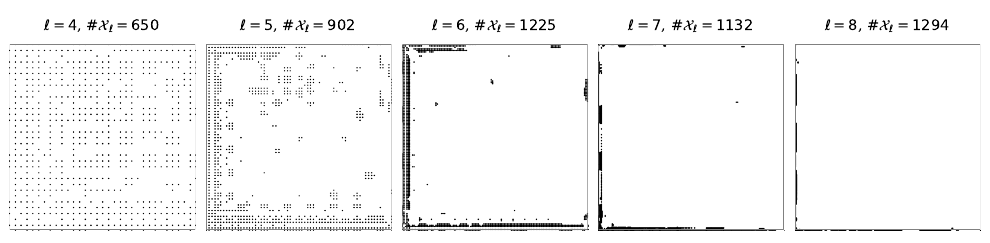}
\caption{Franke's function: selected sets $\frakX_{\ell}$ for the constant threshold $\varepsilon=10^{-4}$, $\ell=4,\ldots,8$.}\label{fig:sets-franke}
\end{figure}

\Cref{fig:sets-franke} shows the selected sets on the finer levels for $\varepsilon=10^{-4}$. While on the coarse levels the selected sites follow the features of Franke's function, on the fine levels they concentrate in a thin layer along the boundary of $\Omega$, where kernel-based interpolants are less accurate. Hence, the remaining sites, roughly $10^3$ per level, are spent on the boundary layer, and the method automatically detects the region where the error of the classical method is concentrated. On the last two levels, only the left and the lower edge remain, where Franke’s function is considerably larger than on the right and the upper edge, $\max|f|\approx0.85$ compared with $0.27$. This is a consequence of the absolute criterion in \cref{alg:adaptive}. Close to the boundary, the error is proportional to the size of the target function.

Finally, \cref{tab:nu} shows the influence of the overlap parameter. A larger $\nu$ improves the classical method considerably. With the constant threshold $\varepsilon=10^{-4}$, the final error and the number of centers of the adaptive method hardly depend on $\nu$, since the error is determined by the threshold. Consequently, the saving decreases with $\nu$, and for $\nu=6$ the classical method reaches the same accuracy with fewer centers. Note that the saving is computed from the final level. Since the error stagnates after about six levels and $\rho_{\ell}$ closely follows the error, a stopping criterion based on $\rho_{\ell}$, which can be evaluated at no extra cost for nested sites, would avoid the last levels and increase the saving. Similarly, the thresholds could be chosen adaptively from the computed values of $\rho_{\ell}$, for instance as $\varepsilon_{\ell} = c \rho_{\ell-1}$. \Cref{thm:aposteriori} remains valid for such choices, since it holds for arbitrary selections of sites. We leave a systematic study to future work.

\begin{table}[ht]
\centering
\small
\begin{tabular}{c cc ccc}
\toprule
& \multicolumn{2}{c}{classical} & \multicolumn{3}{c}{constant, $\varepsilon=10^{-4}$}\\
\cmidrule(lr){2-3}\cmidrule(lr){4-6}
$\nu$ & reduction & rms error & centers & rms error & saving\\
\midrule
2 & 0.27 & $2.5\cdot10^{-6}$ & 6\,844 & $3.2\cdot10^{-5}$ & 3.67\\
3 & 0.23 & $7.1\cdot10^{-7}$ & 5\,512 & $2.8\cdot10^{-5}$ & 2.17\\
4 & 0.20 & $2.6\cdot10^{-7}$ & 5\,312 & $2.5\cdot10^{-5}$ & 1.30\\
6 & 0.18 & $9.2\cdot10^{-8}$ & 5\,182 & $2.4\cdot10^{-5}$ & 0.77\\
\bottomrule
\end{tabular}
\caption{Franke's function with $\phi_{2,1}$, $\mu=1/2$ and $L=8$ for different overlap parameters $\nu$. For the classical method, the mean reduction of the root mean square error per level over the last four levels and the final root mean square error are given. For \cref{alg:adaptive} with the constant threshold $\varepsilon=10^{-4}$, the total number of centers, the final root mean square error and the saving as defined in \cref{tab:thresholds} are given.}\label{tab:nu}
\end{table}

\subsection{The Coupled Regime}\label{subsec:CoupledRegime}

\Cref{thm:main} requires the overlap parameter to be fixed before the refinement parameter is chosen small. To test whether the method deteriorates if $\mu$ and $\nu$ are coupled, we fix $\gamma=\mu\nu\in\{1.5,3\}$ and use $\mu\in\{1/2,1/3,1/4\}$ with $\nu=\gamma/\mu$ for Franke's function. For each pair we use the geometric thresholds with $\vartheta$ equal to the mean error reduction of the classical method, and the constant threshold $\varepsilon=10^{-4}$. The number of levels is chosen such that the finest level has about $6\cdot10^4$ sites.

\Cref{tab:coupled} shows the results. With the geometric thresholds, the adaptive method reproduces the reduction of the classical method for all values of $\mu$, also for $\nu=12$, and reaches a final error within a factor of $1.4$ of the classical one with $56$ to $78$ percent of the centers. With the constant threshold, the final error and $\rho_L$ are about $0.2\varepsilon$ to $0.3\varepsilon$ for all $\mu$, and $\eta_L$ remains of the order of the threshold. Hence, we do not observe any deterioration as $\mu\to0$ with fixed $\gamma$. This suggests that the restriction in \cref{thm:main} is caused by our proof technique, at least for smooth target functions, and strengthens the conjecture in \cref{subsec:discussion}.

\begin{table}[ht]
\centering
\scriptsize
\setlength{\tabcolsep}{3pt}
\begin{tabular}{ccc c cc ccc c ccc}
\toprule
& & & & \multicolumn{2}{c}{classical} & \multicolumn{4}{c}{geometric} & \multicolumn{3}{c}{constant, $\varepsilon=10^{-4}$}\\
\cmidrule(lr){5-6}\cmidrule(lr){7-10}\cmidrule(lr){11-13}
$\gamma$ & $\mu$ & $\nu$ & $L$ & red. & rms error & $\vartheta$ & red. & rms error & centers & rms$/\varepsilon$ & $\rho_L/\varepsilon$ & $\eta_L/\varepsilon$\\
\midrule
1.5 & $1/2$ & 3 & 8 & 0.22 & $1.2\cdot10^{-6}$ & 0.183 & 0.23 & $1.5\cdot10^{-6}$ & 0.56 & 0.27 & 0.27 & 1.03\\
1.5 & $1/3$ & 4.5 & 5 & 0.08 & $9.6\cdot10^{-7}$ & 0.0569 & 0.08 & $1.2\cdot10^{-6}$ & 0.70 & 0.25 & 0.25 & 1.56\\
1.5 & $1/4$ & 6 & 4 & 0.04 & $9.3\cdot10^{-7}$ & 0.0316 & 0.04 & $1.2\cdot10^{-6}$ & 0.67 & 0.20 & 0.20 & 3.58\\
3 & $1/2$ & 6 & 8 & 0.18 & $3.5\cdot10^{-7}$ & 0.152 & 0.17 & $4.8\cdot10^{-7}$ & 0.67 & 0.24 & 0.24 & 1.32\\
3 & $1/3$ & 9 & 5 & 0.06 & $3.0\cdot10^{-7}$ & 0.0417 & 0.05 & $3.7\cdot10^{-7}$ & 0.78 & 0.30 & 0.30 & 1.46\\
3 & $1/4$ & 12 & 4 & 0.03 & $2.9\cdot10^{-7}$ & 0.0214 & 0.03 & $4.1\cdot10^{-7}$ & 0.75 & 0.20 & 0.20 & 1.54\\
\bottomrule
\end{tabular}
\caption{Franke's function with $\phi_{2,1}$ in the coupled regime $\nu=\gamma/\mu$. The number of levels $L$ is chosen such that the finest level has about $6\cdot10^4$ sites. For the classical method and for the geometric thresholds $\varepsilon_\ell=10^{-2}\vartheta^\ell$, where $\vartheta$ is the observed reduction of the classical method, the mean reduction (red.) of the root mean square error per level over the last three levels (two levels for $\mu=1/4$) and the final root mean square error are given. The column centers gives the number of centers of the geometric run relative to the classical method. The last three columns give the final root mean square error, $\rho_L$ and $\eta_L=\|\ce_L\|_{\ell_\infty(\frakD_L)}$ for the constant threshold $\varepsilon=10^{-4}$, divided by $\varepsilon$.}\label{tab:coupled}
\end{table}

\subsection{The L-shaped Domain}

Our second example is the corner singularity
 \begin{align*}
f(\xx)=-r^{2/3}\sin\left(\frac{2t-\pi}{3}\right),
\end{align*}
where we used polar coordinates $ \xx = r(\cos t,\sin t)^\transpose $, with $ t\in[\pi/2,5\pi/2)$, on the L-shaped domain $\Omega=[-1/2,1/2]^2\setminus(0,1/2]^2$, see \cite[Section~6]{Wendland2010}. This function belongs to $H^{\sigma}(\Omega)$ for all $\sigma<5/3$. For $\phi_{2,0}$ the function belongs to the native space, while for $\phi_{2,1}$ we are in the setting of \cref{thm:escape}. The sets of sites are the grids of the previous example, shifted to $[-1/2,1/2]^2$ and restricted to $\Omega$. Hence, $N_\ell=(2^{\ell+1}+1)^2-4^\ell$, i.e., $N_1=21$ and $N_8=197\,633$. The classical method reduces the error by a factor of about $0.30$ per level for both kernels, which is close to $\mu^{5/3}\approx0.31$.

\begin{figure}[ht]
\centering
\includegraphics[width=\textwidth]{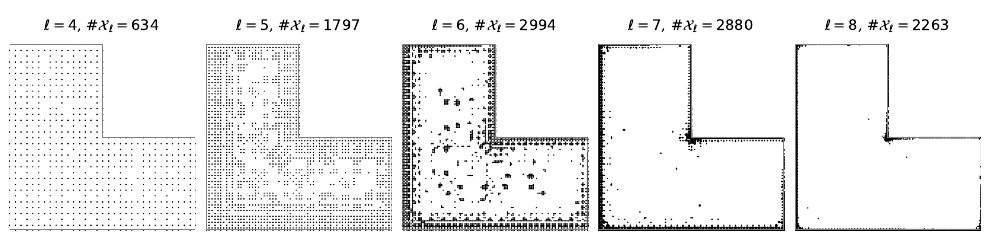}
\caption{L-shaped domain with $\phi_{2,0}$: selected sets $\frakX_{\ell}$ for the constant threshold $\varepsilon=10^{-4}$, $\ell=4,\ldots,8$.}\label{fig:sets-lshape}
\end{figure}

\Cref{tab:lshape} summarizes the results. Again, the error stagnates at about $0.2\varepsilon$ to $0.3\varepsilon$ for constant thresholds. The savings are larger than for Franke's function and reach $5.2$ without and $8.2$ with the check-loop for $\phi_{2,1}$. The check-loop guarantees $\eta_L\leq\tau_L$. The maximum error, however, is attained at the reentrant corner and is essentially the same for all methods. \Cref{fig:sets-lshape} shows that on the finer levels the selected sites concentrate at the reentrant corner and, as for Franke's function, along the boundary of $\Omega$.

\begin{table}[ht]
\centering
\small
\begin{tabular}{ll rcccc}
\toprule
kernel & method & centers & rms error & max error & $\eta_L$ & saving\\
\midrule
$\phi_{2,0}$ & classical & 264\,188 & $5.7\cdot10^{-6}$ & $1.5\cdot10^{-3}$ & -- & --\\
 & constant, $\varepsilon=10^{-3}$ & 2\,074 & $2.2\cdot10^{-4}$ & $1.9\cdot10^{-3}$ & $1.9\cdot10^{-3}$ & 2.07\\
 & constant, $\varepsilon=10^{-4}$ & 10\,812 & $2.9\cdot10^{-5}$ & $1.2\cdot10^{-3}$ & $1.2\cdot10^{-3}$ & 3.89\\
 & check-loop, $\varepsilon=10^{-3}$ & 2\,558 & $1.8\cdot10^{-4}$ & $1.5\cdot10^{-3}$ & $5.0\cdot10^{-4}$ & 2.07\\
 & check-loop, $\varepsilon=10^{-4}$ & 16\,673 & $1.8\cdot10^{-5}$ & $1.5\cdot10^{-3}$ & $5.0\cdot10^{-5}$ & 4.33\\
\midrule
$\phi_{2,1}$ & classical & 264\,188 & $5.8\cdot10^{-6}$ & $2.3\cdot10^{-3}$ & -- & --\\
 & constant, $\varepsilon=10^{-3}$ & 1\,614 & $3.2\cdot10^{-4}$ & $3.3\cdot10^{-3}$ & $3.3\cdot10^{-3}$ & 1.92\\
 & constant, $\varepsilon=10^{-4}$ & 8\,841 & $2.8\cdot10^{-5}$ & $3.3\cdot10^{-3}$ & $3.3\cdot10^{-3}$ & 5.18\\
 & check-loop, $\varepsilon=10^{-3}$ & 1\,914 & $1.7\cdot10^{-4}$ & $2.3\cdot10^{-3}$ & $5.0\cdot10^{-4}$ & 3.20\\
 & check-loop, $\varepsilon=10^{-4}$ & 9\,737 & $1.7\cdot10^{-5}$ & $2.3\cdot10^{-3}$ & $5.0\cdot10^{-5}$ & 8.23\\
\bottomrule
\end{tabular}
\caption{Corner singularity on the L-shaped domain with $\phi_{2,0}$ and $\phi_{2,1}$, $\mu=1/2$, $\nu=3$ and $L=8$: the classical method and \cref{alg:adaptive} with constant thresholds, without and with the check-loop with $\tau_\ell=\varepsilon/2$. Given are the total number of centers, the root mean square and maximum error of $\cf_L$, $\eta_L=\|\ce_L\|_{\ell_\infty(\frakD_L)}$, and the saving as defined in \cref{tab:thresholds}.}\label{tab:lshape}
\end{table}

\subsection{The Composite Function with Distinct Features}\label{subsec:CompositeFunction}

\begin{figure}[ht]
\centering
\includegraphics[width=\textwidth]{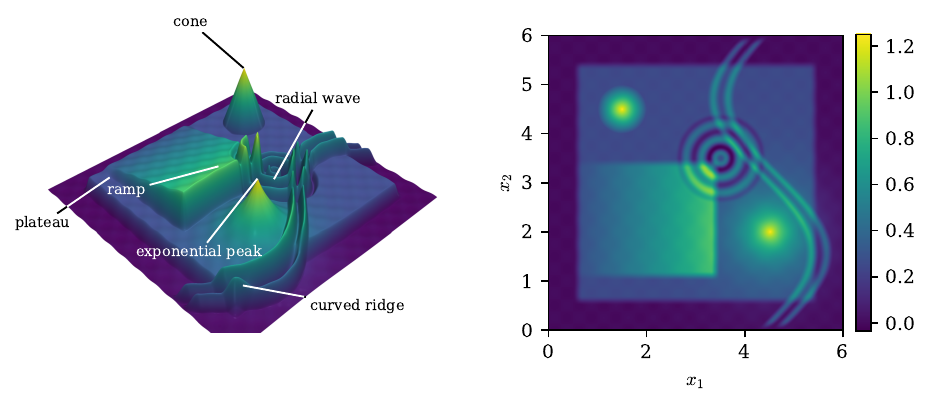}
\caption{The composite function on $\Omega=[0,6]^2$, as a surface (left) and from above (right). It consists of a plateau with smooth edges, a linear ramp with steep but smooth edges, a radial wave, an exponential peak with a kink at its centre, a double ridge along a curve, a cone with kinks at its apex and along its base, and a small oscillation on the whole domain.}
\label{fig:composite-function}
\end{figure}

Our last example is the composite test function of \cite{AvesaniIskeMultererWenzel2026} on $\Omega=[0,6]^2$, see \Cref{fig:composite-function}. It consists of six distinct features, among them a cone, whose gradient jumps along a circle, a ridge along a curve, and several smooth plateaus, whose transitions are only once continuously differentiable. The function belongs to $H^{\sigma}(\Omega)$ only for $\sigma<3/2$, and its singularities are located along curves. We use $\phi_{2,1}$, i.e., the setting of \cref{thm:escape}, and uniform grids with $N_{\ell} = (2^{(\ell+3)} + 1)^2$, i.e., $N_1=289$ and $N_8=4\,198\,401$. The classical method reduces the error per level by factors which approach $\mu^{3/2}\approx0.35$.

\begin{figure}[ht]
\centering
\includegraphics[width=\textwidth]{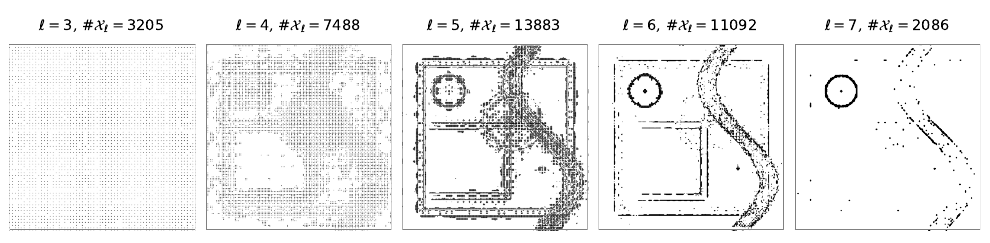}
\caption{Composite function: selected sets $\frakX_{\ell}$ for the constant threshold $\varepsilon=10^{-3}$, for levels  $\ell=3,\ldots,7$.}\label{fig:sets-aimw}
\end{figure}

\Cref{tab:aimw} shows the results per level. On the coarse levels, the adaptive method without check-loop is less accurate than the classical method. Here, the threshold is small compared to the residual, and the update at the dropped sites leads to large residuals, e.g., $\eta_2\approx0.74$ for $\varepsilon=10^{-3}$. On the finer levels, the adaptive method catches up, and the number of selected sites decreases. The check-loop avoids the loss on the coarse levels, since it adds those sites at which the update produced large residuals. For the geometric thresholds, $\vartheta=0.5$ is larger than the classical reduction, such that the adaptive method converges more slowly, compare \cref{tab:thresholds}. \Cref{fig:sets-aimw} shows that the selected sets follow the singularities of the target function: the circle of the cone, the edges of the plateaus, and the ridge.

\begin{table}[ht]
\centering
\scriptsize
\setlength{\tabcolsep}{3pt}
\begin{tabular}{rr c rc rc rcc}
\toprule
& & classical & \multicolumn{2}{c}{constant, $\varepsilon=10^{-3}$} & \multicolumn{2}{c}{constant, $\varepsilon=10^{-4}$} & \multicolumn{3}{c}{check-loop, $\varepsilon=10^{-3}$}\\
\cmidrule(lr){3-3}\cmidrule(lr){4-5}\cmidrule(lr){6-7}\cmidrule(lr){8-10}
$\ell$ & $N_{\ell}$ & rms error & $\#\frakX_{\ell}$ & rms error & $\#\frakX_{\ell}$ & rms error & $\#\frakX_{\ell}$ & passes & rms error\\
\midrule
1 & 289 & $7.9\cdot10^{-2}$ & 289 & $7.9\cdot10^{-2}$ & 289 & $7.9\cdot10^{-2}$ & 289 & 0 & $7.9\cdot10^{-2}$\\
2 & 1\,089 & $5.2\cdot10^{-2}$ & 766 & $6.3\cdot10^{-2}$ & 796 & $6.2\cdot10^{-2}$ & 1\,087 & 2 & $5.2\cdot10^{-2}$\\
3 & 4\,225 & $6.5\cdot10^{-3}$ & 3\,205 & $1.2\cdot10^{-2}$ & 3\,403 & $1.2\cdot10^{-2}$ & 4\,076 & 4 & $6.5\cdot10^{-3}$\\
4 & 16\,641 & $8.0\cdot10^{-4}$ & 7\,488 & $1.7\cdot10^{-3}$ & 11\,826 & $1.4\cdot10^{-3}$ & 10\,533 & 3 & $8.1\cdot10^{-4}$\\
5 & 66\,049 & $1.8\cdot10^{-4}$ & 13\,883 & $4.7\cdot10^{-4}$ & 31\,387 & $2.7\cdot10^{-4}$ & 19\,083 & 5 & $2.3\cdot10^{-4}$\\
6 & 263\,169 & $5.3\cdot10^{-5}$ & 11\,092 & $2.9\cdot10^{-4}$ & 52\,375 & $8.3\cdot10^{-5}$ & 6\,845 & 4 & $1.7\cdot10^{-4}$\\
7 & 1\,050\,625 & $1.7\cdot10^{-5}$ & 2\,086 & $2.8\cdot10^{-4}$ & 59\,232 & $3.8\cdot10^{-5}$ & 4\,732 & 3 & $1.6\cdot10^{-4}$\\
8 & 4\,198\,401 & $6.1\cdot10^{-6}$ & 2\,093 & $2.8\cdot10^{-4}$ & 22\,518 & $2.9\cdot10^{-5}$ & 2\,777 & 6 & $1.6\cdot10^{-4}$\\
\bottomrule
\end{tabular}
\caption{Composite function of \cite{AvesaniIskeMultererWenzel2026} with $\phi_{2,1}$, $\mu=1/2$ and $\nu=3$. For each level $\ell$, the number of sites $N_\ell$, the number of selected sites $\#\frakX_\ell$ and the root mean square error of $\cf_\ell$ are given for the classical method, for \cref{alg:adaptive} with the constant thresholds $\varepsilon=10^{-3}$ and $\varepsilon=10^{-4}$, and for the check-loop with $\varepsilon=10^{-3}$ and $\tau_\ell=5\cdot10^{-4}$. For the check-loop, $\#\frakX_\ell$ is the number of selected sites after the loop, and passes is the number of times the local interpolant was recomputed on level $\ell$.}\label{tab:aimw}
\end{table}

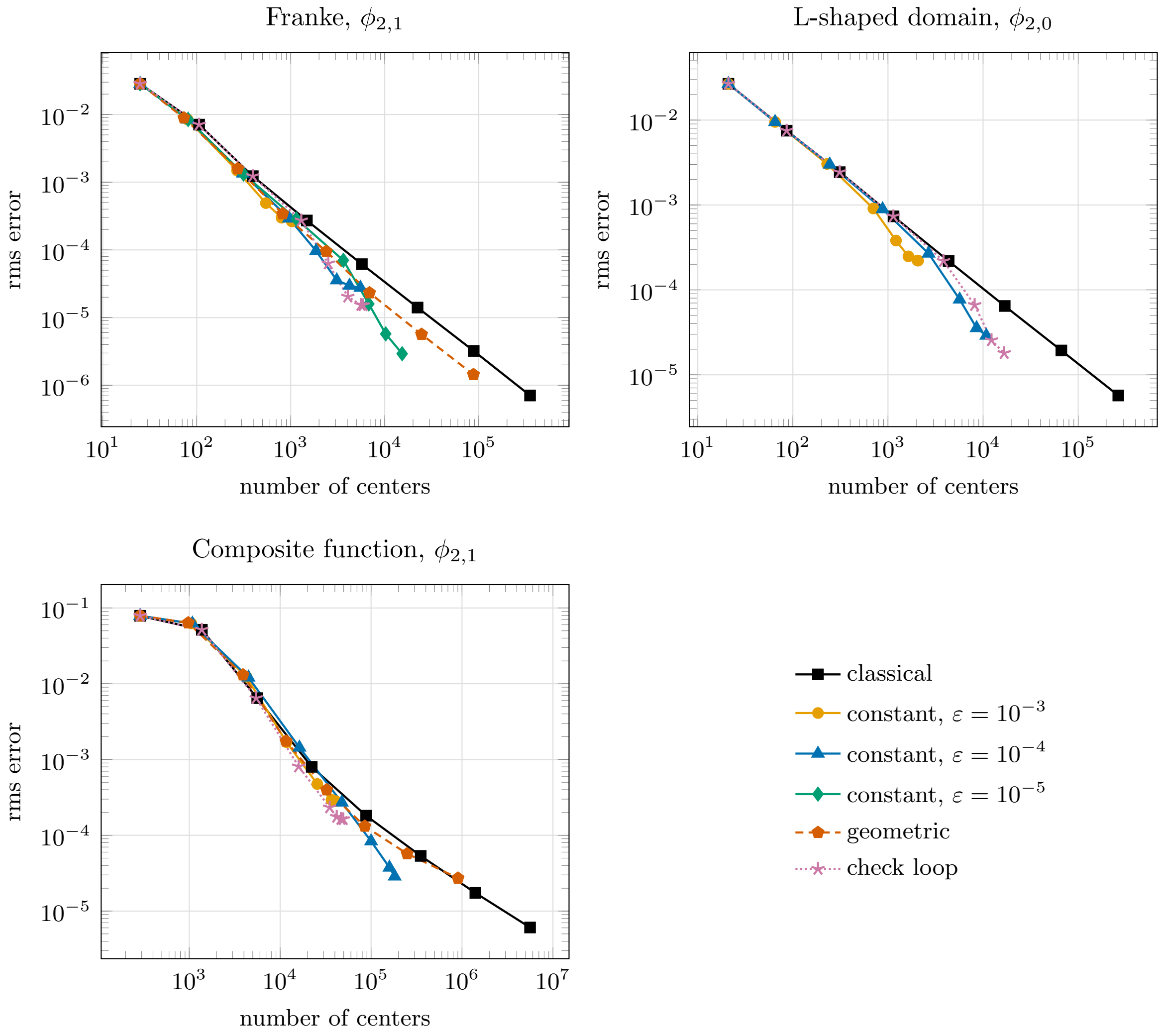
\begin{figure}[ht]
\centering
\providecommand{\figdata}{figures/data}
\definecolor{cClassical}{RGB}{0,0,0}
\definecolor{cEpsA}{RGB}{230,159,0}     
\definecolor{cEpsB}{RGB}{0,114,178}     
\definecolor{cEpsC}{RGB}{0,158,115}     
\definecolor{cGeo}{RGB}{213,94,0}       
\definecolor{cLoop}{RGB}{204,121,167}   
\pgfplotsset{
  amlclassical/.style={cClassical, thick, mark=square*, mark size=1.8pt},
  amlepsA/.style={cEpsA, thick, mark=*, mark size=1.8pt},
  amlepsB/.style={cEpsB, thick, mark=triangle*, mark size=2.3pt},
  amlepsC/.style={cEpsC, thick, mark=diamond*, mark size=2.3pt},
  amlgeo/.style={cGeo, thick, densely dashed, mark=pentagon*, mark size=2.1pt, mark options={solid}},
  amlloop/.style={cLoop, thick, densely dotted, mark=star, mark size=2.6pt, mark options={solid}},
}
\begin{tikzpicture}
\begin{groupplot}[
  group style={group size=2 by 2, horizontal sep=1.6cm, vertical sep=2.1cm},
  width=0.5\textwidth, height=0.42\textwidth,
  xmode=log, ymode=log,
  xlabel={number of centers}, ylabel={rms error},
  grid=major, grid style={gray!25},
  tick label style={font=\footnotesize}, label style={font=\footnotesize},
  title style={font=\small, yshift=-0.5ex},
  table/x=centers, table/y=rms,
]
\nextgroupplot[title={Franke, $\phi_{2,1}$}]
\addplot[amlclassical] table {figures/data/franke_classical.dat};
\addplot[amlepsA]      table {figures/data/franke_c1e-3.dat};
\addplot[amlepsB]      table {figures/data/franke_c1e-4.dat};
\addplot[amlepsC]      table {figures/data/franke_c1e-5.dat};
\addplot[amlgeo]       table {figures/data/franke_geo0.3.dat};
\addplot[amlloop]      table {figures/data/franke_loop.dat};

\nextgroupplot[title={L-shaped domain, $\phi_{2,0}$}]
\addplot[amlclassical] table {figures/data/lshape_classical.dat};
\addplot[amlepsA]      table {figures/data/lshape_c1e-3.dat};
\addplot[amlepsB]      table {figures/data/lshape_c1e-4.dat};
\addplot[amlloop]      table {figures/data/lshape_loop.dat};

\nextgroupplot[title={Composite function, $\phi_{2,1}$}]
\addplot[amlclassical] table {figures/data/comp_classical.dat};
\addplot[amlepsA]      table {figures/data/comp_c1e-3.dat};
\addplot[amlepsB]      table {figures/data/comp_c1e-4.dat};
\addplot[amlgeo]       table {figures/data/comp_geo0.5.dat};
\addplot[amlloop]      table {figures/data/comp_loop.dat};

\nextgroupplot[hide axis, xmin=1, xmax=2, ymin=1, ymax=2,
  legend style={at={(0.5,0.5)}, anchor=center, draw=none, font=\footnotesize,
                cells={anchor=west}, row sep=2pt}]
\addlegendimage{amlclassical}\addlegendentry{classical}
\addlegendimage{amlepsA}\addlegendentry{constant, $\varepsilon=10^{-3}$}
\addlegendimage{amlepsB}\addlegendentry{constant, $\varepsilon=10^{-4}$}
\addlegendimage{amlepsC}\addlegendentry{constant, $\varepsilon=10^{-5}$}
\addlegendimage{amlgeo}\addlegendentry{geometric}
\addlegendimage{amlloop}\addlegendentry{check loop}
\end{groupplot}
\end{tikzpicture}
\caption{Root mean square error against the total number of centers. The check-loop uses $\varepsilon=10^{-4}$ for Franke's function and the L-shaped domain and $\varepsilon=10^{-3}$ for the composite function; the geometric thresholds use $\vartheta=0.3$ for Franke's function and $\vartheta=0.5$ for the composite function.}\label{fig:error-centers}
\end{figure}

\Cref{fig:error-centers} summarizes the three examples by plotting the root mean square error against the total number of centers.

\subsection{Discussion}

Our numerical results confirm the qualitative predictions of the theory and show that the rate of \cref{thm:main} is pessimistic. For constant thresholds, the error stagnates at about $0.2\varepsilon$ to $0.3\varepsilon$ in all examples, as predicted by \cref{cor:thresholds}. The quantity $\rho_{\ell}$, which enters the a posteriori bound of \cref{thm:aposteriori}, follows the error closely, while the residual at the dropped sites in the maximum norm can be much larger than the threshold. For geometric thresholds, however, we observe the rate $\max(\alpha,\vartheta)$ instead of the rate $\max(r,\vartheta)$ guaranteed by \cref{cor:thresholds}, i.e., the adaptive method attains the reduction of the classical method if $\vartheta$ does not exceed the classical reduction. This supports \cref{conj:rate}, and the best choice is $\vartheta$ close to the classical reduction. The adaptive method reaches a given accuracy with up to eight times fewer centers. The savings are largest for the target functions with singularities and for small thresholds. The selected sets concentrate at the singularities and along the boundary of the domain and hence indicate the locations of singularities and boundary layers, which may be of interest for feature detection. Finally, we did not observe any deterioration of \cref{alg:adaptive} in the coupled regime $\nu=\gamma/\mu$, which is not covered by \cref{thm:main}, but only by \cref{cor:checkloop-rate} for the variant with check-loop.

\section{Summary and Outlook}\label{sec:summary}

We have proven the convergence of the truly adaptive kernel-based multilevel method, in which the local interpolants are computed only on those sites at which the residual of the preceding level exceeds a threshold the user can choose. The key ingredients are the $\ell_2$-version of the sampling inequality, the uniform localization of the Lagrange functions on arbitrary subsets of a quasi-uniform set of sites, and the smoothness of the residual on the coarser scale. As a by-product, we obtained an a posteriori error bound for arbitrary point selections, which can be evaluated at negligible cost for nested sites, and a convergence result for a variant of the method with an additional check-loop, which converges with the rate of the classical multilevel method. Furthermore, the analysis carries over to penalized least-squares approximation and to target functions which are rougher than the native space of the kernel.

We see several directions for future research. First, our a priori bound requires the overlap parameter to be fixed before the refinement parameter is chosen small, excluding the coupled setting $\nu\geq\gamma/\mu$. The key obstructions are discussed in detail in \cref{subsec:discussion}: For the full sets $X_{\ell}$, a localization which is uniform in $\nu$ can be expected from \cite{HangelbroekNarcowichRiegerWard2018}. For subsets with holes larger than $\delta_{\ell}$, however, the Lebesgue constants grow with $\nu$, since the Lagrange functions extrapolate beyond the boundary of the holes. Hence, an extension to the coupled setting cannot rely on bounds which are uniform with respect to all subsets, but has to use properties of the residual close to the boundary of the selected sets. In the numerical tests in \cref{subsec:CoupledRegime} we did not see any deterioration in the coupled setting. Second, the rate $\sqrt\mu$ is most likely not optimal. In the numerical examples, the adaptive method attains the reduction of the classical method, $\mu^s+\nu^{-s}$ up to a constant close to one, if the threshold decreases at the same rate. Finally, the combination of the adaptive multilevel method with the tensor product multilevel method \cite{BuettnerKempfWendland2026,KempfWendland2023} and its direction-adaptive variants should be investigated.

\section*{Acknowledgements}
The author used the AI assistant Claude (Anthropic) during the preparation of this work, for revising the text, and for writing parts of the code for the numerical experiments. The author reviewed, verified and edited all content, and takes full responsibility for the content of this article.

\bibliographystyle{abbrv}
\bibliography{references}

\end{document}